\documentclass[12pt,reqno]{article}
\usepackage[english]{babel}
\usepackage[T1]{fontenc}
\usepackage{graphics}
\usepackage{graphicx}
\usepackage{amsfonts}
\usepackage[all]{xy}
\usepackage{amsmath}
\newtheorem{theorem}{Theorem}[section]
\newtheorem{proposition}{Proposition}[section]
\newtheorem{corollary}{Corollary}[section]
\newtheorem{lemma}{Lemma}[section]

\makeatletter \@addtoreset{equation}{section}

\newtheorem{remark}{Remark}[section]
\newcommand{\C}{\mathbb C}
\newcommand{\R}{\mathbb R}
\newcommand{\Z}{\mathbb Z}

\begin{document}

\title
{ The Hua operators on Homogeneous Line Bundles \\ over Bounded Symmetric Domains of Tube Type.}

\author {Abdelhamid Boussejra
\thanks{e-mail: a.boussejra@gmail.com}  \\
Department of
Mathematics,  Faculty of Sciences\\
 University Ibn Tofail, Kenitra
, Morocco. }
\date{}
\maketitle
\begin{abstract}
Let $\mathcal{D}=G/K$ be a bounded symmetric domain of tube type. We show that the image of the Poisson transform on the degenerate principal series representation of $G$  attached to the Shilov boundary of $\mathcal{D}$ is characterized by a $K$-covariant differential operator on a homogeneous line bundle over $\mathcal{D}$. As a consequence of our result we get  the eigenvalues of the Casimir operator for Poisson transforms on homogeneous line bundles over $G/K$. This extends a result of Imemura and all \cite {I} on symmetric domains of classical type to all symmetric domains. Also we compute  a class of Hua type integrals generalizing an earlier result of Faraut and Koranyi\cite{F}.
 \end{abstract}
\section{Introduction and main result}
 Let $G/K$ be an irreducible Hermitian symmetric space of tube type with Shilov boundary $G/P_{\Xi}$. In \cite{Sh} Shimeno proved that the Poisson transform maps the space $B(G/P_{\Xi},L_{\lambda})$ of huperfunction-valued sections of a degenerate spherical series representation attached to $G/P_{\Xi}$ bijectively onto an eigenspace of the Hua operator $\mathcal{H}$ on $G/K$, under certain conditions on the complex parameter $\lambda$.\\
 The aim of this paper is to generalize the above result to  homogeneous line bundles over $G/K$.\\
To state our result in rough form let us fix some notations refering to section 2 for more details.\\
For $\nu\in {\Z}$ let $E_{\nu}$ denote the homogeneous line bundle on $G/K$ associated to the one dimensional representation $\tau_{\nu}$ of $K$. Let $\mathbb{D}_{\nu}(G/K)$ be the algebra of $G$-invariant differential operators on $E_{\nu}$. Shimeno \cite{Sh1} proved that the Poisson transform  is a $G$-isomorphism from the space $\mathcal{B}(G/P,L_{\mu,\nu})$ of hyperfunction-valued sections of principal series representations attached to the Furstenberg boundary $G/P$ onto the solutions space $\mathcal{A}(G/K,\mathcal{M}_{\mu,\nu})$ of the system of differential equations on $E_{\nu}$
$$
\mathcal{M}_{\mu,\nu}:\quad (D-\chi_{\mu,\nu}(D))F=0\quad \forall D\in \mathbb{D}_{\nu}(G/K),
$$  
under certain conditions on $\nu$ and $\mu\in \mathfrak{a}^{\ast}_{c}$. In above $\chi_{\mu,\nu}$ is a certain character of the algebra $\mathbb{D}_{\nu}(G/K)$.\\ 
Let $r$ and $m$ denote respectively the rank of $G/K$ and the multiplicity of the short restricted roots.
Let $P_{\Xi}$ be a maximal standard parabolic subgroup of $G$ with the Langlands decomposition $P_{\Xi}=M_{\Xi}A_{\Xi}N_{\Xi}$ such that $A_{\Xi}$ is of real dimension one. Let $\xi$ be the one dimensional representation of $P_{\Xi}$ defined by
 $$
 \xi(m_{1}man)=\tau_{\nu}(m)a^{\rho_{\Xi}-\lambda\rho_{0}},\quad m_{1}\in M_{\Xi,s}, m\in M, a\in A_{\Xi}, n\in N_{\Xi}.
$$
In above $M_{\Xi,s}$ is the semisimple part of $M_{\Xi}$, $\rho_{0}$ and $\rho_{\Xi}$ are  linear forms on $\mathfrak{a}_{\Xi}={\R}X_{0}$ (the Lie algebra of $A_{\Xi}$) defined by $\rho_{0}(X_{0})=r$ and $\rho_{\Xi}=\eta \rho_{0}$ with $\eta=\frac{m}{2}(r-1)+1$.\\ 

Let $L_{\xi}$ be the  homogeneous line bundle over the Shilov boundary $G/P_{\Xi}$ associated to $\xi$.\\
For $f$ in  $B(G/P_{\Xi},L_{\lambda,\nu})$  the space of hyperfunction-valued sections of the homogeneous line bundle $L_{\xi}$ we define its Poisson transform by
 $$
 [P_{\lambda,\nu}f](g)=\int_{K} f(gk)\tau_{\nu}(k)dk,
 $$
 where $dk$ is the normalized Haar measure of $K$.\\
 
 The  degenerate series representation $B(G/P_{\Xi},L_{\lambda,\nu})$ attached to $G/P_{\Xi}$ is a $G$-sub-module of a principal series representation $B(G/P,L_{\mu_{\lambda},\nu})$ for  $\mu_{\lambda}=\lambda\rho_{0}-\rho_{\Xi}+\rho$, and the image $P_{\lambda,\nu}(B(G/P_{\Xi},L_{{\lambda},\nu}))$ is a $G$-submodule of the solution space $\mathcal{A}(G/K,\mathcal{M}_{\mu_{\lambda},\nu})$.\\
Therefore it is natural to pose the problem of characterizing this image by differential operators on $E_{\nu}$.

In the trivial cases the origin of this problem   goes back to L. H. Hua \cite{H} who showed in the case of the classical Cartan domain of $n\times n$ matrices that for $f$ in $B(G/P_{\Xi})$ $P_{\rho_{\Xi}}f$   is annihilated by $n^{2}$-second order differential operators. Since then many authors considered the problem of constructing differential operators characterizing the image $P_{\rho_{\Xi}}(B(G/P_{\Xi},)$, see \cite{M}, \cite{J} in the case of the Siegel upper half plane. In  \cite{Jk} Johnson and Koranyi constructed second order differential operator $\mathcal{H}$-called then after Hua operator- and showed that  $\mathcal{H}$ characterizes the image $P_{\rho_{\Xi}}(B(G/P_{\Xi},)$.
Lassalle \cite{L} reproved their result introducing the operator $\mathcal{H}_{\mathfrak{q}}$, cutting down the number of equations.\\
In this paper we will show that the image of $B(G/P_{\Xi},L_{\lambda,\nu})$ under $P_{\lambda,\nu}$ is  characterized by a \\$K$- covariant differential operator on $E_{\nu}$ .\\

Let $\mathfrak{g}$ and $\mathfrak{k}$ be the Lie algebras of $G$ and $K$ respectively and $\mathfrak{g}=\mathfrak{k}\oplus \mathfrak{p}$ be the Cartan decomposition of $G$ with Cartan involution $\theta$. We denote  the complexifications of $\mathfrak{g}, \mathfrak{k}$, $ \mathfrak{p}$ by $\mathfrak{g}_{c}, \mathfrak{k}_{c}$, $ \mathfrak{p}_{c}$   respectively.\\
 The center $\mathfrak{z}$ of $\mathfrak{k}$ is of dimension one and there exists $Z_{0}\in \mathfrak{z}$ such that $adZ_{0}$ define a complex structure on $\mathfrak{p}_{c}$. Let
$$
\mathfrak{g}_{c}=\mathfrak{p}_{-}\oplus \mathfrak{k}_{c}\oplus \mathfrak{p}_{+},
$$
be the corresponding eigenspace decomposition of $ \mathfrak{p}_{c}$.\\
Let ${E_{i}}$ be a basis of $\mathfrak{p}_{+}$ and ${E_{i}^{\ast}}$ be the dual basis of $\mathfrak{p}_{-}$ with respect to the Killing form of $\mathfrak{g}_{c}$.\\
We consider the element of $\mathcal{U}(\mathfrak{g}_{c})\otimes \mathfrak{k}_{c}$ -called here the  Hua operator- defined by
\begin{equation}\label{Hua}
 \mathcal{H}=\sum_{i,j}E_{i}E_{j}^{\ast}\otimes [E_{j},E_{i}^{\ast}].
\end{equation}
The operator $\mathcal{H}$ is a homogeneous differential operator from the space of $C^{\infty}$-sections of $E_{\nu}$ to the space of $C^{\infty}$-sections of the homogeneous vector bundle on $G/K$ associated to the representation $\tau_{\nu}\otimes Ad_{K}\vert_{{\mathfrak{k}_{c}}}$, which does not depend on the choice of basis. (see section 4)\\
The pair $(K,K\cap M_{\Xi})$ is a compact symmetric pair. Let $\widetilde{\tau}$ be the corresponding involution of $\mathfrak{k}$,  and let $\mathfrak{k}=\mathfrak{l}\oplus \mathfrak{q}$ be the decomposition of $\mathfrak{k}$ into eigenspaces of  $\widetilde{\tau}$. Here $\mathfrak{l}$ is the Lie algebra of $K\cap M_{\Xi}$\\
Let $\mathcal {H}_{q}$ be the element of $\mathcal{U}(\mathfrak{g}_{c})\otimes \mathfrak{q}_{c}$ defined by
\begin{equation}\label{Hua1}
\mathcal{H}_{q}=\sum_{i,j}E_{i}E_{j}^{\ast}\otimes \mathrm{p}([E_{j},E_{i}^{\ast}]),
\end{equation}
where $ \mathrm{p}$ denotes the orthogonal projection of $\mathfrak{k}_{c}$ onto $\mathfrak{q}_{c}$.\\

The main result of this paper can be stated as follows
\begin{theorem}
Let $\lambda$ be a complex number and let $\nu\in {\Z}$ such that 
\begin{equation}
-\lambda-\frac{m}{2}(-r+2+j)\notin \{1,2...\} \quad \textit{for} \quad j=0,1,
\end{equation}
\begin{equation}
-\lambda+\eta-\mid\nu\mid\notin 2{\Z}^{+}+2.
\end{equation}
Then the Poisson transform $P_{\lambda,\nu}$ is a $G$-isomorphism from the space $B(G/P_{\Xi},L_{\lambda,\nu})$ onto the space of $C^{\infty}$- sections of $E_{\nu}$ satisfying the system of differential equations
\begin{eqnarray}
\mathcal{H}_{q}F=\frac{(\lambda^{2}-(\eta-\nu)^{2})}{4p}F.(-iZ_{0}).
\end{eqnarray}
\end{theorem}
In the above  $p=2\eta$ is the genus of the bounded symmetric domain $G/K$.\\

In the case $\tau_{\nu}$ is the trivial representation, Theorem 1.1 has been established  by Lassalle \cite{L} for $\lambda=\rho_{\Xi}$ and generalizes to generic $\lambda$ by Shimeno \cite{Sh}.\\  

To prove our result we first show that every solution of the Hua system is a joint eigenfunction of the algebra $\mathbb{D}_{\nu}(G/K)$ (Theorem 6.1).\\To do so, let $\mathfrak{h}\subset \mathfrak{q}$ be a Cartan subalgebra of the symmetric pair $(\mathfrak{k},\mathfrak{l})$ and let $\mathfrak{s}$ be the orthogonal complement of $\mathfrak{h}$ in $ \mathfrak{q}$ with respect to the Killing form $B$. Write 
$$
\mathcal{H}_{\mathfrak{q}}=\mathcal{H}_{\mathfrak{h}} +\mathcal{H}_{\mathfrak{s}} 
$$

with $\mathcal{H}_{\mathfrak{h}}\in \mathcal{U}(\mathfrak{g}_{c})\otimes \mathfrak{h}_{c}$ and $\mathcal{H}_{\mathfrak{s}}\in\mathcal{U}(\mathfrak{g}_{c})\otimes \mathfrak{s}_{c}$.\\

Then we prove (Theorem 6.2) that if F is a $\tau_{-\nu}$-spherical function on $G$  satisfying 
$$
\mathcal{H}_{\mathfrak{h}}F=\frac{(\lambda^{2}-(\eta-\nu)^{2})}{4p}F.(-iZ_{0}),  
$$
then the function
$$
\phi(t_{1},...,t_{r})=\prod_{j=1}^{r}(\cosh t_{j})^{\nu}F(t_{1},...,t_{r}),
$$
satisfies the following system of differential equations
$$
\frac{\partial^{2}\phi}{\partial t^{2}_{k}}+2\coth2t_{k}\frac{\partial\phi}{\partial t_{k}}-2\nu\tanh t_{k}\frac{\partial\phi}{\partial t_{k}}
+\frac{m}{2}\sum^{r}_{j=1}\frac{1}{(\sinh^{2}t_{j}-\sinh^{2}t_{k})}(\sinh 2t_{j}\frac{\partial\phi}{\partial t_{j}}-\sinh 2t_{k}\frac{\partial\phi}{\partial t_{k}})
$$
$$
=\frac{(\lambda^{2}-(\eta-\nu)^{2})}{4}\phi,
$$
for all $k=1,...,r$.\\
Then by using a result of Yan \cite{Y} on generalized hypergeometric functions in several variables we deduce that $F$ is given (up to a constant) in terms of the generalized Gauss-hypergeometric function $_{2}F_{1}^{(m)}(b,c,d;x_{1},....,x_{r})$ associated to the parameter $m$, see \cite{Y}. Namely
$$
F(t_{1},....,t_{r})=\prod _{j=1}^{r}(1-\tanh^{2}t_{j})^{\frac{\lambda+\eta}{2}}\quad_{2}F_{1}^{(m)}(\frac{\lambda+\eta-\nu}{2},
\frac{\lambda+\eta+\nu}{2};\eta,\tanh^{2}t_{1},...,\tanh^{2}t_{r}).
$$
To finish the proof of our main result we show that under the conditions (1.3) and (1.4) the induced  equations of the subsystem $\mathcal{H}_{\mathfrak{s}}F=0$  for boundary values on the maximal boundary $G/P$ characterize the space $\mathcal{B}(G/P_{\Xi},L_{\lambda,\xi})$.\\
\subsection{Consequences.}
For the symmetric domains of classical type, the eigenvalues of the Casimir operator for Poisson transforms on $B(G/P_{\Xi},L_{\lambda,\nu})$ have been computed by Imemura and all \cite{I}. Using our result we can compute them for all symmetric domains without using their classification as in \cite{I}. Thus, of course the exceptional domains are included. \\ More precisely, let $L$ be the Casimir operator acting on $C^{\infty}$-sections of $E_{\nu}$. 
\begin{corollary}
Let $\lambda\in {\C}, \nu\in{\Z}$ and let $F=P_{\lambda,\nu}f$ with $f\in \mathcal{B}(G/P_{\Xi},L_{\lambda,\xi})$. Then we have
\begin{equation}
LF=\frac{\lambda^2-(\eta-\nu)^2}{4r}F. 
\end{equation}
\end{corollary}
\textbf{Proof}.
The Laplacian $L$ is the projection of the Hua operator $\mathcal{H}$ onto the center $\mathfrak{z}_{c}$ of $\mathfrak{k}_{c}$. Thus $\mathcal{H}_{\mathfrak{z}}=iL\otimes Z_{0}^{\ast}$.\\
Since $\mathfrak{z}\subset \mathfrak{h}$, then 
\begin{equation}
\mathcal{H}_{\mathfrak{z}}F=\frac{\lambda^2-(\eta-\nu)^2)}{4p}F (-iZ_{0}),
\end{equation}
by Theorem 1.1.\\
Let $B$ be the Killing form of $\mathfrak{g}$. Then we have 
$$
B(\mathcal{H}_{\mathfrak{z}}F,Z_{0})=iLF,
$$
from which we deduce that $LF= -\frac{\lambda^2-(\eta-\nu)^2)}{4p}F B(Z_{0},Z_{0})$. Since $B(Z_{0},Z_{0})=-2\eta r$, the result follows.\\ 
\begin{remark}
The above corollary agrees with [\cite{K2},Theorem 5.2 with $q=\eta$ and  $s=\eta-\lambda$ ] where the eigenvalues of the Casimir operator have been computed by Koranyi using a different method.
\end{remark}    
Let $h(z,w)$ be the canonical Jordan polynomial associated to the bounded symmetric domain $G/K$ and let $S$ denote its Shilov boundary.\\
Then by using trivialization of the homogeneous line bundles $E_{\nu}$ and $L_{\xi}$ we can rewrite the Poisson transform   
$P_{\lambda,\nu}$ as
\begin{equation}\label{Poisson integral} 
P_{\lambda,\nu}f(z)=\int_{S}[\frac{h(z,z)}{\mid h(z,u)\mid^{2}}]^{\frac{\lambda+\eta-\nu}{2}}h(z,u)^{-\nu}f(u)du,
\end{equation}
see Proposition 5.1, for more details.\\

Letting $f=1$ in (\ref{Poisson integral})  we deduce from  the method of the proof of Theorem 1.1, an explicit expression of a class of  Hua type integrals. Namely, we have
\begin{corollary}
let $\lambda\in{\C}$ and let $\nu\in{\Z}$. Then we have
$$
\int_{S}[\frac{h(z,z)}{\mid h(z,u)\mid^2}]^{\frac{\lambda+\eta-\nu}{2}}h(z,u)^{-\nu}du=h(z,z)^{\frac{\lambda+\eta-\nu}{2}}
$$
$$
\quad_{2}F_{1}^{(m)}(\frac{\lambda+\eta-\nu}{2},\frac{\lambda+\eta+\nu}{2};\eta,\tanh^{2}t_{1},...,\tanh^{2}t_{r}),
$$
$z=Ad(k)\sum\limits_{j=1}^{r}\tanh t_{j}E_{\gamma_{j}}$.
\end{corollary}
The above formula has been established for $\nu=0$ by Faraut and Koranyi in \cite{F}.\\

In the case of the trivial line bundle on $SU(n,n)/S(U(n)\times U(n))$, we computed a more general class of the following Hua type integrals. More precisely, we have 
$$
\int_{U(n)}[\frac{h(\tanh t I ,\tanh tI)}{\mid h(\tanh tI,u)\mid}]^{\frac{\lambda+\eta}{2}}\phi_{\mathbf{m}}(u)du=\frac{n!}{d_{\mathbf{m}}}\det (\phi_{\lambda,\mid m_{i}-i+j\mid}(t))_{i,j},
$$
where $\mathbf{m}=(m_{1},...,m_{n})\in {\Z}^{n}$ such that $m_{1}\ge m_{2}\ge...\ge m_{n}$ and $d_{\mathbf{m}}=\Pi_{i<j}(1+\frac{m_{i}-m_{j}}{j-i})$, with  
$$
\phi_{\lambda,k}(t)=(1-\tanh^{2}t)^{\frac{\lambda+n}{2}}\frac{(\frac{\lambda+n}{2})_{k}}{(1)_{k}}\tanh^{k}t\quad _{2}F_{1}(\frac{\lambda+n}{2},\frac{\lambda+n}{2}+k,1+k;\tanh^{2}t).      
$$   
 In above $_{2}F_{1}(a,b,c;x)$ is the classical Gauss hypergeometric function, $(a)_{k}$ the Pochammer symbol and  $\phi_{\mathbf{m}}$ a zonal spherical function related to Schur functions, see \cite{B} for more details.\\

Before ending this section, we should mention that 
Professor Koufany informed me that he and Professor Zhang had also obtained a characterization of  Poisson integrals on homogeneous line bundles over tube type symmetric domains see \cite{K1}, See also \cite{So}.\\ Namely, they showed that the image $P_{\lambda,\nu}(B(G/P_{\Xi},L_{\lambda,\nu}))$ can be  characterized by the operator $\mathcal{H}$ defined in(\ref{Hua}). The system (\ref{Hua}) has of course $dim_{\R}K$ differential equations.\\ Regarding to the scalar case, this system has too many equations.\\
Our result shows that, as in the scalar case, a subsystem of $dim_{\R} S=dim_{\C}G/K$ equations is sufficient to characterize Poisson integrals on homogeneous line bundles over $S$ .\\ 

The organization of this paper is as follows. After a preliminaries on Hermitian symmetric spaces we review in section 3 the results of \cite{Sh} on the Poisson transform on homogeneous line bundles on the Furstenberg boundary. In section 4 we define the Poisson transform on degenerate principal series representation attached to the Shilov boundary and introduce the Hua operator $\mathcal{H}$ on $E_{\nu}$. Using a trivialization of the space of $C^{\infty}$-sections of $E_{\nu}$ we give a realization of $\mathcal{H}$ on the Harish-Chandra realization of $G/K$ (Proposition 4.1). Section 5, 6 and 7 are devoted to the proof of our main result.\\

Acknowledgement: I would like to thank Professor Koranyi for sending to me his recent work \cite{K2} , where the necessity of the Hua equations is also proved and  for several valuable comments on a first version of this paper.

\section{Preliminaries and notations}
In this section we recall some structural results on Hermitian symmetric space, see \cite{H} for more details.\\
For a real Lie algebra $\mathfrak{b}$ we shall denote by $\mathfrak{b}_{c}$ its complexification.\\
Let $G/K$ be an irreducible hermitian symmetric space of noncompact type with rank $r$. The group $G$ is the connected component of its isometry group and $K$ is the isotropy subgroup of $G$.\\
Let $\mathfrak{g}=\mathfrak{k}\oplus \mathfrak{p}$ be  the Cartan decomposition of the Lie algebra $\mathfrak{g}$ of $G$ with respect to a Cartan involution $\theta$.Thus $\mathfrak{k}$ has one dimensional center $\mathfrak{z}$, $\mathfrak{k}_{s}=[\mathfrak{k},\mathfrak{k}]\neq \mathfrak{k}$ and $\mathfrak{k}=\mathfrak{z}\oplus \mathfrak{k}_{s}$.\\ 
Let $Z_{0}\in \mathfrak{z}$ such that $(ad Z_{0})^{2}=-1$ on $\mathfrak{p}_{c}$. Let $\mathfrak{p}_{+}$ (respectively $\mathfrak{p}_{-}$) be the $i$ (respectively $-i$) -eigenspace of $ad Z_{0}$ in $\mathfrak{g}_{c}$.\\
Then $ \mathfrak{p}_{+}$ and $\mathfrak{p}_{-}$ are Abelian subalgebras. Moreover 
$
[\mathfrak{p}_{+},\mathfrak{p}_{-}]=\mathfrak{k}_{c}.
$
We thus have the Harish-Chandra decomposition of $\mathfrak{g}_{\C}$
$$
\mathfrak{g}_{\C}=\mathfrak{p}_{+}+ \mathfrak{k}_{\C} +\mathfrak{p}_{-}.
$$
Let $G_{\C}$ be the simply connected Lie group with Lie algebra $g_{c}$. 
Let $P_{+}, P_{-}$ and $K_{\C}$  denote the analytic subgroup of $G_{\C}$ corresponding to the Lie subalgebras $\mathfrak{p}_{+}, \mathfrak{p}_{-}$ and $ \mathfrak{k}_{\C}$ respectively. Then $P_{+}K_{c} P_{-}$ is an open dense subset of $G_{c}$ containing $G$. For $z\in\mathfrak{p}_{+}$ and $g\in G$ we denote by $U(g:z)$ the $K_{\C}$ component of $g\exp z$. That is
\begin{eqnarray}
g\exp z=\exp(g.z) U(g:z)p_{-}(g).
\end{eqnarray} 
Under the above action the $G$-orbit $\mathcal{D}=G.0$ of $z=0\in \mathfrak{p}_{+}$ is a bounded domain in $\mathfrak{p}_{+}$ and $K$ is the isotropy subgroup of $0$. This is the Harish-Chandra realization of the Hermitian symmetric space $G/K$.\\
Recall that the $K_{\C}$-component $U(g:z)$ is the canonical automorphy factor and that it satisfies the multiplier identity
\begin{equation}
U(g_{1}g_{2}:z)=U(g_{1}:g_{2}.z)U( g_{2}:z).
\end{equation}
\subsection{The Roots.}
Let $\mathfrak{t}$ be a Cartan subalgebra of $\mathfrak{k}$ (and hence also of $\mathfrak{g}$). Let $\Delta$ denote the root system of  $(\mathfrak{g}_{\C},\mathfrak{t}_{c})$. \\

For $\alpha\in \Delta$ let $\mathfrak{g}_{\alpha}$ denote the root space for $\alpha$. A root $\alpha$ is said to be compact (resp. noncompact) if the root space $\mathfrak{g}_{\alpha}$ is contained in $\mathfrak{k}_{\C}$ (resp.$\mathfrak{p}_{c}$).
Let $B$ denote the Killing form of $\mathfrak{g}_{\C}$. For each   root $\alpha$ we can choose  root vectors $\tilde{H_{\alpha}}\in \mathfrak{t}_{c}$, $\tilde{E_{\alpha}}\in \mathfrak{g}_{\alpha}$, and $\tilde{E_{-\alpha}}\in \mathfrak{g}_{-\alpha}$ such that
$$
\alpha(H)=B(H,\tilde{H}_{\alpha}),\quad \forall H\in \mathfrak{t}_{c},
$$
$$
[\tilde{E_{\alpha}},\tilde{E_{-\alpha}}]= \tilde{H_{\alpha}},  \quad \tau \tilde{E_{\alpha}}=-\tilde{E_{-\alpha}},
$$
with $B(\tilde{E_{\alpha}},\tilde{E_{-\alpha}})=1$.\\
In above  $\tau$ denotes the conjugation in $\mathfrak{g}_{\C}$ with respect to the real form $\mathfrak{k}+i\mathfrak{p}$.\\

For $\alpha,\beta$ in $\Delta$, we set $<\alpha,\beta>=B(\tilde{H_{\alpha}},\tilde{H_{\beta}})$. Then the length $\mid \alpha\mid$ of a root $\alpha$ is defined  by $\mid \alpha\mid=\sqrt{<\alpha,\alpha>}$. \\
Put $c_{\alpha}=\frac{\sqrt{2}}{\mid\alpha\mid}$. Let  $H_{\alpha}= c_{\alpha}^{2}\tilde{H_{\alpha}}$ and   $E_{\alpha}=c_{\alpha}\tilde{E_{\alpha}}$.
Then the root vectors $H_{\alpha},E_{\alpha}$ and $E_{-\alpha}$ satisfy
$$
[E_{\alpha},E_{-\alpha}]=H_{\alpha},\quad [E_{\alpha},E_{-\alpha}]=H_{\alpha},
$$
$$
B(E_{\alpha},E_{-\alpha})=\frac{2}{\mid\alpha\mid^{2}}.
$$ 
Moreover  $\alpha(H_{\alpha})=2$.\\

We choose an ordering on $\Delta$ such that a noncompact root is positive if and only if $\mathfrak{g}_{\alpha}\subset \mathfrak{p}_{+}$.
We denote by $\Phi^{+}$ the set of positive noncompact roots.Then we have
$$
\mathfrak{p}_{+}=\sum_{\alpha\in \Phi^{+}}\mathfrak{g}_{\alpha},
$$
$$
\mathfrak{p}_{-}=\sum_{\alpha\in \Phi^{+}}\mathfrak{g}_{-\alpha}.
$$
For each $\alpha\in \Phi^{+}$ we set 
$$
X_{\alpha}=E_{\alpha}+E_{-\alpha},\quad  Y_{\alpha}=i(E_{\alpha}-E_{-\alpha}).
$$
Two roots $\alpha$ and $\beta$ are said strongly orthogonal if neither $\alpha+\beta$ nor $\alpha-\beta$ are roots.\\
Let $\Gamma=\{\gamma_{1},....,\gamma_{r}\}$ be a maximal set of strongly orthogonal  noncompact roots, such that $\gamma_{j}$ is the highest element of $\Phi^{+}$ strongly orthogonal to $\gamma_{j+1},....,\gamma_{r}$, for $j=r,....,1$. Then $\mathfrak{a}=\sum_{j=1}^{r}{\R}X_{\gamma_{j}}$
is a maximal Abelian subspace of $\mathfrak{p}$.\\
Let $\mathfrak{h}$ denote  the Abelian subalgebra generated by the elements $iH_{\gamma_{j}}, j=1,...,r$ and let  $\mathfrak{h}^{\perp}$ be its orthogonal in $\mathfrak{t}$ with respect to the Killing form $B$. Then $\mathfrak{h}^{\perp}=\{H\in \mathfrak{t};\gamma_{j}(H)=0, j=1,...,r\}$.\\

For $\alpha,\beta\in \Delta$ denote $\alpha\sim \beta$ if and only if $\alpha\mid_{\mathfrak{h}}=\beta\mid_{\mathfrak{h}}$. Let
$$
\Phi^{+}_{ij}=\{\alpha\in \Phi^{+};\alpha\sim\frac{\gamma_{i}+\gamma_{j}}{2}\},
$$
for $1\leq i<j\leq r$.\\

Let $C$ denote the set of compact roots in $\phi$. Define
$$
C_{0}=\{\alpha\in C;\alpha\sim 0\},
$$
and
$$
C_{ij}=\{\alpha\in C;\alpha\sim\frac{\gamma_{j}-\gamma_{i}}{2}\}, \hspace*{.8cm}\quad \mbox{for}   1\leq i<j\leq r.
$$
Then it is known that $\Delta^{+}$ is the disjoint union of the sets $ \Gamma, C_{0},C_{i,j}$, $\Phi_{ij}^{+}$ and $\Phi^{+}$ is the disjoint union of the sets $\Gamma$ and $\Phi^{+}_{ij}$.\\

Let $\alpha\in \Phi^{+}_{ij}$. Define  $\widetilde{\alpha}\in \mathfrak{t}^{\star}_{c}$
by $\widetilde{\alpha}=\gamma_{i}+\gamma_{j}-\alpha $. Then $\widetilde{\alpha}\in \Phi_{ij}^{+}$ see \cite{L}.\\
Now we recall a result from \cite {B-V} which will be  helpful later.

\begin{proposition} Let $k$ be a fixed integer, $1\leq k\leq r$.\\
(i) For $\gamma_{j}\in \Gamma$ we have $ \gamma_{j}(H_{\gamma_{k}})=2\delta_{jk}$.\\
(ii) If $\alpha\in \Phi^{+}_{jk}$, then $\alpha (H_{\gamma_{k}})=1$.\\
(iii) In all other cases $\alpha(H_{\gamma_{k}})=0$.\\
(iv)Let $\alpha\in \Phi^{+}_{jk}$. Then $<\alpha,\alpha>=<\gamma_{k},\gamma_{k}>$ if $\alpha\neq\widetilde{\alpha}$ and $<\alpha,\alpha>=\frac{1}{2}<\gamma_{k},\gamma_{k}>$ if $\alpha=\widetilde{\alpha}$.
\end{proposition}

Let $c$ be the Cayley transform of $\mathfrak{g_{c}}$ given by
$$
c=\exp\frac{\pi}{4}(\overline{E}_{0}-E_{0}),
$$
where $E_{0}=\sum_{k=1}^{r}E_{\gamma_{k}}$.
Then $Ad c(H_{\gamma_{j}})=X_{\gamma_{j}}$, and  $Ad c(\mathfrak{h})=i\mathfrak{a}$.

It is well known that $Ad c^{4}=1$ and  $Ad c^{2}$   is an automorphism of $\mathfrak{k}$. Let $\mathfrak{l}$ (respectively $\mathfrak{q}$) be the $+1$ (respectively $-1$)eigenspace of $Ad c^{2}$ in $\mathfrak{k}$. Then we have $\mathfrak{z}$ is a subset of $\mathfrak{q}$ . Moreover, let  $\mathfrak{q}_{s}=\mathfrak{q}\cap \mathfrak{k}_{s}$ then $\mathfrak{k}_{s}=\mathfrak{q}_{s}+\mathfrak{l}$ and $\mathfrak{q}= \mathfrak{q}_{s}+\mathfrak{z}$.\\

Put $\beta_{i}=\gamma_{j}\circ (Ad c^{-1}|a)$, for $j=1,...,r$. Then the set of restricted roots $\Sigma$ of the pair $(\mathfrak{g},\mathfrak{a})$ is given by
$$
\Sigma=\{\pm \beta_{j}\quad (1\le j\le r), \frac{\pm \beta_{j}\pm \beta_{k}}{2} \quad (1\le j\neq k\le r)\}.
$$
Let $\Sigma^{+}$ be the set of positive roots in $\Sigma$. Then
$$
\Sigma^{+}=\{ \beta_{j}\quad (1\le j\le r), \frac{ \beta_{j}\pm \beta_{k}}{2} \quad (1\le k<j\le r)\}.
$$
The Weyl group $W$ of $\Sigma$ acts as the group of all permutations and sign changes of the set $\{\beta_{1},....,\beta_{r}\}$, so it is isomorphic to the semi-direct product of $({\Z}/2{\Z})^{r}$ and the symmetric group.\\
We set
$$
\alpha_{j}=\frac{\beta_{r-j+1}-\beta_{r-j}}{2},\quad (1\le j\le r-1)\quad \alpha_{r}=\beta_{1}.
$$
Then  $\Gamma=\{\alpha_{1},...,\alpha_{r}\}$ is the set of simple roots in $\Sigma^{+}$. Let $\{H_{1},...,H_{r}\}$ denote the basis of $\mathfrak{a}$ which is dual to $\{\alpha_{1},...,\alpha_{r}\}$.\\
For $\alpha\in \Sigma$ let $\mathfrak{g}^{\alpha}$ denote the corresponding root space and let $m_{\alpha}$ be the multiplicity of $\alpha$.\\
The multiplicities of the roots $\frac{\pm \beta_{j}\pm \beta_{k}}{2}$ and $\pm \beta_{j}$ are $m$ and $1$ respectively.\\
As usual set $\rho=\frac{1}{2}\sum_{\alpha\in \Sigma^{+}}m_{\alpha}\alpha$, $\mathfrak{n}^{+}=\sum_{\alpha\in \Sigma^{+}}\mathfrak{g}^{\alpha}$ and $\mathfrak{n}^{-}=\theta(\mathfrak{n})$.\\

Let $A,N^{+}$ and $N^{-}$ be the analytic subgroups of $G$ corresponding to $\mathfrak{a},\mathfrak{n}^{+}$ and $\mathfrak{n}^{-}$ respectively.The group $G$ has the Iwasawa decomposition $G=KAN^{+}$. Let $M$ be the centralizer of $\mathfrak{a}$ in $K$. Then $P=MAN^{+}$ is a minimal parabolic subgroup of $G$.\\

 Let $\Xi=\Gamma\setminus\{\alpha_{r}\}$
and let $P_{\Xi}$ be the corresponding standard parabolic subgroup of $G$ with the Langlands decomposition $P_{\Xi}=M_{\Xi}A_{\Xi}N_{\Xi}$ such that $A_{\Xi}\subset A$. Then $P_{\Xi}$ is a maximal standard parabolic subgroup of $G$ and the space $G/P_{\Xi}$ is the Shilov boundary of $X$.\\
If $\mathfrak{a}_{\Xi}$ denotes the Lie algebra of $A_{\Xi}$, then
$$
\mathfrak{a}_{\Xi}=\{H\in \mathfrak{a}; \gamma(H)=0, \forall\gamma\in\Xi\}.
$$
Moreover $\mathfrak{a}_{\Xi}={\R}X_{0}$ where $X_{0}=\sum_{j=1}^{r} X_{\gamma_{j}}$.\\
On $\mathfrak{a}_{\Xi}$ we define the linear form $\rho_{0}$ by $\rho_{0}(X_{0})=r$.
Let $\rho_{\Xi}$ be the restriction of $\rho$ to $\mathfrak{a}_{\Xi}$. Then
$$
 \rho_{\Xi}=(m\frac{r-1}{2}+1)\rho_{0}.
$$
The algebras $\mathfrak{n}^{\pm}$ decomposes as $\mathfrak{n}^{\pm}=\mathfrak{n}_{\Xi}^{\pm}+\mathfrak{n}({\Xi})^{\pm}$, where
$$
\mathfrak{n}_{\Xi}^{\pm}=\sum_{j\le k}\mathfrak{g}^{\pm(\frac{\beta_{j}+\beta_{k}}{2})}, \quad \mathfrak{n}({\Xi})^{\pm}=\sum_{j>k}\mathfrak{g}^{\pm(\frac{\beta_{j}-\beta_{k}}{2})},
$$
and we have $\mathfrak{p}_{\Xi}=m+a+n^{+}+\mathfrak{n}({\Xi})^{-}$, $\mathfrak{p}_{\Xi}$ been the Lie algebra of $P_{\Xi}$.\\

\section{Eigensections of invariant differential operators}

We review the main result of \cite{Sh1} on the image of the Poisson transform on the principal series representation attached to the Furstenberg boundary G/P.\\
Recall that    the genus $p$ of the bounded domain $\mathcal{D}$, is given by
$$
p=m(r-1)+2.
$$ 
Then the length $\mid \gamma_{j}\mid$  of the roots $\gamma_{j}$ is such that $\mid \gamma_{j}\mid=\frac{1}{\sqrt{p}}$.\\
Let $Z=\frac{p}{n}Z_{0}$ where $n=dim_{\C} \mathfrak{p}_{+}$.\\
Let $K_{s}$ be the analytic subgroup of $K$ with Lie algebra $\mathfrak{k}_{s}$. For $\nu\in {\Z}$ define $\tau_{\nu}:K\rightarrow {\C}^{\times}$ by 
$$
\tau_{\nu}(k)=1 \quad \textit{if} \quad k\in K_{s},
$$ 
and 
$$
\tau_{\nu}(\exp(tZ))=e^{-i\nu t}\quad \textit{for} \quad t\in{\R}.
$$
Then $\tau_{\nu}$ determines a one dimensional representation of $K$ and all one dimensional representations of $K$ have this form, see \cite{Sc} and the remark below.\\
\begin{remark}
Since 
$$
B(Z,Z_{0})=-2p,
$$
it follows that  $Z$ is the same as the element
$$
\frac{1}{r}\sum^{r}_{j=1}\frac{2i\tilde{H}_{\gamma_{j}}}{<\gamma_{j},\gamma_{j}>},
$$
in \cite{Sc}.   
\end{remark}

Let $E_{\nu}$ be the homogeneous line bundle  on $G/K$ associated to $\tau_{\nu}$.
The space of $C^{\infty}$-sections of $E_{\nu}$ can be identified with the space 
$$
C^{\infty}(G/K,\tau_{\nu})=\{f\in C^{\infty}(G);f(gk)=\tau_{\nu}(k)^{-1}f(g)\quad \textit{for all}\quad g\in G, k\in K\}.
$$
The group $G$ acts on $C^{\infty}(G/K,\tau_{\nu})$ by the left regular representation $\pi(g)f(x)=f(g^{-1}x)$.\\ 
Let $\mathbb{D}_{\nu}(G/K)$ be the set of all left-invariant differential operators on $G$ that map $C^{\infty}(G/K,\tau_{\nu})$ into itself. Then, accordingly to Shimeno result \cite{Sh1} $\mathbb{D}_{\nu}(G/K)$ is isomorphic, via the Harish-Chandra isomorphism  $\gamma_{\nu}$, to $\mathcal{U}(\mathfrak{a})^{W}$ the set of Weyl group invariant elements in $\mathcal{U}(\mathfrak{a})$.\\
For $\mu\in\mathfrak{a}^{\ast}_{c}$ and $\nu\in {\Z}$ we define an algebra homomorphism $\chi_{\mu,\nu}$ of $\mathbb{D}_{\nu}(G/K)$ by
$$
 \chi_{\mu,\nu}(D)=\gamma_{\nu}(D)(\mu).
$$
For $\mu\in \mathfrak{a}^{\ast}_{c}$ we denote by $\mathcal{B}(G/P,L_{\mu,\nu})$ the space of hyperfunction-valued sections of the homogeneous line bundle on $G/P$ associated to the character $\sigma_{\mu,\nu}$ of $P$ given by
$$
\sigma_{\mu,\nu}(man)=a^{\rho-\mu}\tau_{\nu}(m) \quad m\in M, a\in A, n\in N^{+}.
$$
For $f\in \mathcal{B}(G/P,L_{\mu,\nu})$ we define the Poisson transform $\mathcal{P}_{\mu,\nu}$ by
$$
\mathcal{P}_{\mu,\nu}f(g)=\int_{K}f(gk)\tau_{\nu}(k)dk.
$$

A straightforward computation shows that
$$
\mathcal{P}_{\mu,\nu}f(g)=\int_{K}e^{-(\mu+\rho)H(g^{-1}k)}f(k)\tau_{\nu}(\kappa( g^{-1}k))dk,
$$
where $\kappa :G\rightarrow K$ and $H:G\rightarrow \mathfrak{a}$ are the projections defined by $g\in \kappa(g)e^{H(g)}N^{+}$.\\

Let $\mathcal{B}(G/K,\tau_{\nu})$ be the space of hyperfunction-valued sections of the homogeneous line bundle $E_{\nu}$ on $G/K$. We denote by $\mathcal{A}(G/K,\mathcal{M}_{\mu,\nu})$  the space of all real analytic functions in $\mathcal{B}(G/K,\tau_{\nu})$ which satisfy the system of differential equations
$$
\mathcal{M}_{\mu,\nu}: DF=\gamma_{\nu}(D)(\mu)F, \quad D\in \mathbb{D}_{\nu}(G/K).
$$
Let $e_{\nu}(\mu)^{-1}$ be the denominator of the $c$-function associated to $\mathcal{P}_{\mu,\lambda}$. That is
$$
e_{\nu}(\mu)^{-1}=\prod_{1\le j<k\le r}\Gamma (\frac{1}{2}(m+\mu_{j}+\mu_{k}))\Gamma (\frac{1}{2}(m+\mu_{k}-\mu_{j}))
$$
$$
\times \prod_{1\le j\le r}\Gamma (\frac{1}{2}(1+\mu_{j}+\nu))\Gamma (\frac{1}{2}(1+\mu_{j}-\nu)).
$$

\begin{theorem}\cite{Sh}
Let $\mu\in \mathfrak{a}_{c}^{\ast}$ and $\nu\in {\Z}$ satisfying the conditions
$$
-2\frac{<\mu,\alpha>}{<\alpha,\alpha>}\notin \{1,2,...\} \quad \mbox{for all} \quad\alpha\in \Sigma^{+}
$$
and
$$
 e_{\nu}(\mu)\neq 0,
$$
then the Poisson transform $\mathcal{P}_{\mu,\nu}$ is a $G$-isomorphism from $\mathcal{B}(G/P,L_{\mu,\nu})$ onto $\mathcal{A}(G/K,\mathcal{M}_{\mu,\nu})$.\\The inverse of  $\mathcal{P}_{\mu,\nu}$ is given by the boundary value map up to a non-zero constant multiple.
\end{theorem}
Recall that a ${\C}$-valued function on $G$ is called $\tau_{-\nu}$-spherical function if it satisfies
\begin{equation}\label{spherical} 
f(k_{1}gk_{2})=\tau_{-\nu}(k_{2})f(g)\tau_{-\nu}(k_{1})
\end{equation}
A $\tau_{-\nu}$-spherical function $\Phi$ on $G$ will be called an elementary spherical function of type $\tau_{-\nu}$ if it satisfies
$$
D\Phi=\chi_{\mu,\nu}\Phi,\quad \forall D\in \mathbf{D}_{\nu}(G/K)
$$
$$
\Phi(e)=1
$$
According to \cite{Sh}, for $\mu\in \mathfrak{a}^{\ast}_{c}$ and $\nu\in {\Z}$, there exists a unique elementary spherical function of type $\tau_{-\nu}$, given by 
\begin{equation}\label{elementary spherical} 
\Phi_{\mu,\nu}(g)=\int_{K}e^{-(\mu+\rho)H(g^{-1}k)}\tau_{\nu}(k^{-1}\kappa( g^{-1}k))dk.
\end{equation}

\section{The Poisson transform and the Hua operator}

\subsection{The Poisson transform on a homogeneous line bundle}
 In this subsection  we define the Poisson transform  on degenerate principal series representation  attached to the Shilov boundary $G/P_{\Xi}$ of $G/K$.\\
Let $M_{\Xi,s}$ be the analytic subgroup of $M_{\Xi}$ with Lie algebra $[\mathfrak{m}_{\Xi},\mathfrak{m}_{\Xi}]$. From here on we suppose that $\tau_{\nu}$ is such that $\tau_{\nu}\mid_{K\cap M_{\Xi,s}}=1$.\\
For $\lambda\in {\C}$ and $\nu\in {\Z}$ let $\xi_{\lambda,\nu}$ denote  the one dimensional representation of $P_{\Xi}$ defined by 
$$
\xi_{\lambda,\nu}(m_{1}man)=a^{\rho_{\Xi}-\lambda\rho_{0}}\tau_{\nu}(m), \textit{for all} \quad m_{1}\in M_{\Xi,s}, m\in M, a\in A_{\Xi}, n\in N_{\Xi}. 
$$
Let $B(G/P_{\Xi},L_{\lambda,\nu})$ be the space of hyperfunction-valued sections of the line bundle on $G/P_{\Xi}$ associated to the character  $\xi_{\lambda,\nu}$. 
The Poisson transform $P_{\lambda,\nu}$ of an element $f\in B(G/P_{\Xi},L_{\lambda,\nu})$ is defined by
\begin{eqnarray}
P_{\lambda,\nu}f(g)=\int_{K}f(gk)\tau_{\nu}(k)dk.
\end{eqnarray}

By the generalized Iwasawa decomposition $G=KP_{\Xi}$, the restriction map from $G$ to $K$ gives an isomorphism from $B(G/P_{\Xi},L_{\lambda,\nu})$ onto the space $B(K/K_{\Xi},L_{\nu})$ of hyperfunction-valued sections of the homogeneous line bundle  on $K/K_{\Xi}$ associated to the representation $\tau_{\nu}$. Here $K_{\Xi}=K\cap M_{\Xi}$.\\
For $\lambda\in {\C}$ define the following ${\C}$-linear form  $\mu_{\lambda}$ on $ \mathfrak{a}^{\ast}_{c}$ by
$$
\mu_{\lambda}(H)=(\lambda\rho_{0}-\rho_{\Xi})(H_{\Xi})+\rho(H),
$$
where $H_{\Xi}$ is the $\mathfrak{a}_{\Xi}$-component of $H$ with respect to the orthogonal decomposition $\mathfrak{a}=\mathfrak{a}_{\Xi}\oplus \mathfrak{a}({\Xi})$. 
We have
\begin{eqnarray}
B(G/P_{\Xi},L_{\lambda,\nu})\subset B(G/P,L_{\mu_{\lambda},\nu}).
\end{eqnarray}
From (4.2) we deduce 
\begin{eqnarray}
P_{\lambda,\nu}(B(G/P_{\Xi},L_{\lambda,\nu}))\subset \mathcal{A}(G/K,\mathcal{M}_{\mu_{\lambda},\nu}).
\end{eqnarray}

A straightforward computation shows that the Poisson transform of $f\in B(K/K_{\Xi},L_{\nu})$ is given by
 
$$
P_{\lambda,\nu}f(g)=\int_{K}e^{-(\lambda+\eta)\rho_{0}(H_{\Xi}(g^{-1}k))}f(k)\tau_{\nu}(\kappa(g^{-1}k))dk.
$$

The space  $B(K/K_{\Xi},L_{\nu})$ can be identified to the space  of all hyperfunctions $f$ on $K$  such that
$$
f(km)=\tau_{\nu}^{-1}(m)f(k) \qquad \forall\quad m\in K_{\Xi}.
$$
Let $\Lambda$ be the map  from $B(K/K_{\Xi}),L_{\xi})$ into the space $B(K/K_{\Xi})$ of all hyperfunctions on the Shilov boundary $K/K_{\Xi}$ defined by 
$$
\Lambda f(k)=\tau_{\nu}(k)f(k).
$$
Then $\Lambda$ is a $K$-isomorphism.\\
Using $\Lambda$  the Poisson transform (denoted again by $P_{\lambda,\nu}$) of an element $f\in B(K/K_{\Xi})$ is given by
\begin{eqnarray}
P_{\lambda,\nu}f(g)=\int_{K}e^{-(\lambda+\eta)\rho_{0}H_{\Xi}(g^{-1}k)}f(k)\tau_{\nu}(k^{-1}\kappa(g^{-1}k))dk, 
\end{eqnarray}

\subsection{The Hua operator.}
If $(\theta,V)$ is a finite dimensional representation of the compact group $K$, we denote by $C^{\infty}(G/K,\theta)$ the space of $C^{\infty}$-sections of the homogeneous vector bundle on $G/K$ associated to $\theta$.\\
Let ${E_{i}}$ be a basis of $\mathfrak{p}_{+}$ and $E^{\ast}_{i}$ be the dual basis of $\mathfrak{p}_{-}$ with respect to the Killing form $B$.\\
Let $\mathcal{U}(\mathfrak{g}_{\C})$ denote the universal enveloping algebra of $\mathfrak{g}_{\C}$. We consider the element of $\mathcal{U}(\mathfrak{g}_{\C})\otimes \mathfrak{k}_{\C}$    defined by
$$
 \mathcal{H}=\sum_{i,j}E_{i}E_{j}^{\ast}\otimes [E_{j},E_{i}^{\ast}].
$$
Then $\mathcal{H}$ defines a homogeneous differential operator  from the space $C^{\infty}(G/K,\tau_{\nu})$ to the space 
$C^{\infty}(G/K,\tau_{\nu}\otimes Ad_{K}\mid_{\mathfrak{k}_c})$, which does not depend on the choice of the basis.\\
Let $V$ be a linear subspace of $\mathfrak{k}$ and let $V_{c}$ be its complexification. We denote by $p$ the orthogonal projection from
$\mathfrak{k}_{\C}$ onto $V_{c}$. We extend $p$ on $\mathcal{U}(\mathfrak{g}_{\C})\otimes \mathfrak{k}_{\C}$ by setting
$$
p(U\otimes X)=U\otimes p(X) \quad (U\in \mathcal{U}(\mathfrak{g}_{\C}), X\in \mathfrak{k}_{\C}).
$$
We put $\mathcal{H}_{V}=p(\mathcal{H})$.\\
If $v_{j}$ is a basis of $V_{c}$ and $v_{j}^{\ast}$ is the dual basis with respect to $B$. Then
$$
\mathcal{H}_{V}=\sum_{k}U_{k}\otimes v^{\ast}_{k},
$$
where $U_{k}$ is the element of $\mathcal{U}(\mathfrak{g}_{\C})$ given by
$$
U_{k}=\sum_{j}[v_{k},E_{j}]E^{\ast}_{j},
$$
see \cite{L} for more details.\\

Let $\Lambda_{\nu}$ be the operator defined on $C^{\infty}(G/K,\tau_{\nu})$ by
\begin{eqnarray}
\Lambda_{\nu}F(z)=\tau_{\nu}(U(g:0))F(g),\quad z=g.0.
\end{eqnarray}
Then $\Lambda_{\nu}$ is an isomorphism from $C^{\infty}(G/K,\tau_{\nu})$ onto $C^{\infty}(\mathcal{D})$.
Notice that  
$$
\tau_{-\nu}(U(g:z))=[J(g,z)]^{\frac{\nu}{p}},
$$
where $J(g,z)$ stands for the Jacobian of the transformation $g$.\\

Now we define an action $T_{\nu}$ of $G$ on $\mathcal{D}$ as follows:\\
For each $g\in G$ define $T_{\nu}(g)$ such that the following diagram
$$
\xymatrix{
    C^{\infty}(G/K,\tau_{\nu}) \ar[r]^{\pi(g)} \ar[d]_{\Lambda_{\nu}}  & C^{\infty}(G/K,\tau_{\nu}) \ar[d]^{\Lambda_{\nu}} \\
    C^{\infty}(\mathcal{D}) \ar[r]_{T_{\nu}(g)} & C^{\infty}(\mathcal{D})
 }
$$ 
is commutative.\\
The following result can be proved by direct computations
\begin{lemma}
i) Let $F\in C^{\infty}(\mathcal{D})$. For any $g\in G$  we have 
$$
T_{\nu}(g)F(z)=\tau_{-\nu}(U(g^{-1}:z))F(g^{-1}z)
$$
ii) The operator $\Lambda_{\nu}$ is a $G$-intertwining operator from $C^{\infty}(G/K,\tau_{\nu})$ onto $C^{\infty}(\mathcal{D})$.  
\end{lemma}

Using the above $G$-intertwining operator the Hua operator may be viewed as acting  on $C^{\infty}(\mathcal{D})$. The new operator which we denote by $\mathcal{H}_{\nu}$ will be given below.\\
For $f\in C^{\infty}(G/K,\tau_{\nu}\otimes Ad)$, we define the function $\Lambda_{\tau_{\nu}\otimes Ad}f:G/K 
\rightarrow {\C}\otimes \mathfrak{k}_{\C}$ by 
$$
\Lambda_{\tau_{\nu}\otimes Ad}f(z)=(\tau_{\nu}\otimes Ad)(U(g:0))f(g),\quad z=g.0.
$$
We define the Hua operator $\mathcal{H}_{\nu}$ on $C^{\infty}(\mathcal{D})$ such that the following diagram

$$
\xymatrix{
    C^{\infty}(G/K,\tau_{\nu}) \ar[r]^{\mathcal{H}} \ar[d]_{\Lambda_{\nu}}  & C^{\infty}(G/K,\tau_{\nu}\otimes Ad_{K}\mid_{\mathfrak{k}_c}) \ar[d]^{\Lambda_{\nu}\otimes Ad} \\
    C^{\infty}(\mathcal{D}) \ar[r]_{\mathcal{H}_{\nu}} & C^{\infty}(\mathcal{D},{\C}\otimes\mathfrak{k}_{\C})
 }
$$    
is commutative.

\begin{proposition}
Let $F\in C^{\infty}(\mathcal{D})$. Then we have
\begin{equation}
\mathcal{H}_{\nu}F(z)=\tau_{\nu}(U(g:0))\sum_{i,j}[Ad(U(g:0)^{-1})E_{i}Ad(U(g:0)^{-1})E^{\ast}_{j}]\\
\Lambda_{\nu}^{-1}F(g)\otimes [E_{j},E^{\ast}_{i}],
\end{equation}
with $z=g.0$.
\end{proposition}

Proof. The proof follows by direct computations.\\

The Hua operator $\mathcal{H}_{\nu}$ has the following invariance property
\begin{proposition}
For any $h\in G$ and $F\in C^{\infty}(\mathcal{D})$ we have
$$
T_{\nu}(h)(\mathcal{H}_{\nu})F(z)=Ad(U(h^{-1}:z))\mathcal{H}_{\nu}( T_{\nu}(h)F)(z).
$$

\end{proposition}

Proof. Let $g\in G$ such that $g.0=z$. We have
\begin{equation*}\begin{split}
\mathcal{H}_{\nu}F(h^{-1}.z)&=\tau_{\nu}(U(h^{-1}g:0))\\
&\times\sum_{i,j}[Ad(U(h^{-1}g:0)^{-1})E_{i}Ad(U(h^{-1}g:0)^{-1})E^{\ast}_{j}]\Lambda_{\nu}^{-1}F(h^{-1}g)\otimes [E_{j},E^{\ast}_{i}],
\end{split}\end{equation*}
by Proposition 4.1.\\
Use the  identity (2.2) on the factor of automorphy $U(g:z)$ to write  
$$
U(h^{-1}g:0)=U(h^{-1},g.0)U(g:0),
$$
and the identity 
$$
\sum_{i,j}Ad(k^{-1})E_{i}Ad(k^{-1})E^{\ast}_{j}f(g)\otimes[E_{j},E^{\ast}_{i}]=
\sum_{i,j}E_{i}E^{\ast}_{j}f(g)\otimes[Ad(k)E_{j},Ad(k)E^{\ast}_{i}],\quad \forall k\in K_{c},
$$
 to get
\begin{equation}\begin{split}
\mathcal{H}_{\nu}F(h^{-1}.(z))&=\tau_{\nu}(U(h^{-1}:z))Ad (U(h^{-1}:z))\tau_{\nu}(U(g:0))\\
&\times\sum_{i,j}Ad(U(g:0)^{-1})E_{i}Ad(U(g:0)^{-1})E^{\ast}_{j}
\Lambda_{\nu}^{-1}F(h^{-1}g)\otimes[E_{j},E^{\ast}_{i}].
\end{split}\end{equation}
Next by using  the fact that $\Lambda_{\nu}^{-1}F(h^{-1}g)=\Lambda_{\nu}^{-1}(T_{\nu}(h)F)(g)$ as well as the formula (4.6) we obtain
$$
\tau_{-\nu}(U(h^{-1}:z))\mathcal{H}_{\nu}F(h^{-1}.z)=Ad(U(h^{-1}:z))\mathcal{H}_{\nu}((T_{\nu}(h)F)(z),
$$
and the proposition follows.\\

\section{Necessity of the conditions}
In this section we will establish that Poisson integrals are eigenfunctions of the Hua operator.\\

\begin{lemma}
Let $F\in C^{\infty}(G/K,\tau)$. Then 
$$
\mathcal{H}F=\frac{(\lambda^{2}-(\eta-\nu)^{2})}{4p}F.(-iZ_{0}).
$$ 
if and only if\\ 
the function $\tilde{F}=\Lambda_{\nu}F$ satisfies $\mathcal{H}_{\nu}\tilde{F} =\frac{(\lambda^{2}-(\eta-\nu)^{2})}{4p}\tilde{F}.(-i Z_{0})$.
\end{lemma}

Proof. It follows by direct computations.\\

Let $du$ be the normalized $K$-invariant measure on the Shilov boundary $S$. By using the operator $\Lambda_{\nu}$ we can rewrite the Poisson transform (4.4)  from $ B(S)$ to $C^{\infty}(\mathcal{D})$  
$$
P_{\lambda,\nu}f(z)=\int_{S}P_{\lambda,\nu}(z,u)f(u)du,
$$
where
$$
P_{\lambda,\nu}(z,u)=e^{-(\lambda+\eta)\rho_{0}(H_{\Xi}(g^{-1}k))}\tau_{\nu}(U(g:0))\tau_{\nu}(k^{-1}\kappa(g^{-1}k)),
$$
with $z=g.o$ and $u=k.E_{0}$.\\

Next we introduce   a $K$-invariant polynomial $h(z)$ on $\mathfrak{p}_{+}$ whose restriction on $\sum_{j=1}^{r}{\R}E_{\gamma_{j}}$ is given by
$$
h(\sum_{j=1}^{r}a_{j}E_{\gamma_{j}})=\Pi_{j=1}^{r}(1-a_{j}^{2}).
$$
Let $h(z,w)$ denote  its polarization, see \cite{F} for more details. 

\begin{proposition}

i)The Poisson kernels for the homogeneous line bundles are given by
$$
P_{\lambda,\nu)}(z,v)=[\frac{h(z,z)}{\mid h(z,v)\mid^{2}}]^{\frac{\lambda+\eta-\nu}{2}}h(z,v)^{-\nu}.
$$
\begin{eqnarray}
ii)P_{\lambda,\nu}(g.z,g.u)=P_{\lambda,\nu}(z,u)J_{g}(z)^{-\frac{\nu}{2p}}J_{g}(u)^{-\frac{\lambda+\eta-\nu}{2p}}
\overline{J_{g}(u)}^{-\frac{\lambda+\eta+\nu}{2p}},
\end{eqnarray}
for every $g\in G$.
\end{proposition}

Proof. i) Let 
\begin{eqnarray}
\Psi_{\lambda,\nu}(z)=e^{-(\lambda+\eta)\rho_{0}(H_{\Xi}(g^{-1}))}\tau_{\nu}(U(g:0))\tau_{\nu}(\kappa(g^{-1})\quad z=g.0
\end{eqnarray}
Observe that the right hand side of (5.2) is right $K$-invariant, hence $\Psi_{\lambda,\nu}(z)$ is well defined  on $\mathcal{D}$. 
Define $\mu_{\lambda}\in \mathfrak{a}^{\ast}_{c}$ by  $\mu_{\lambda}=\lambda\rho_{0}+\rho_{\Xi}-\rho$. Then 
$$
e^{(\lambda\rho_{0}+\rho_{\Xi})H_{\Xi}(g^{-1})}=e^{(\mu_{\lambda}+\rho)H(g^{-1})}.
$$

Next recall that if $g=ne^{A(g)}\kappa_{1}(g)$  with respect to the decomposition $G=NAK$, then 
$$
A(g)=-H(g^{-1}), \kappa_{1}(g)=(\kappa(g^{-1}))^{-1}.
$$
Henceforth 
\begin{eqnarray}
\Psi_{\lambda,\nu}(g)=e^{(\mu_{\lambda}+\rho)A(g)}\tau_{\nu}(U(g:0))\tau_{-\nu}(\kappa_{1}(g)).
\end{eqnarray}
But the right hand-side of (5.3) is nothing but the generalized Harish-Chandra $c$-function $e_{\lambda,\nu}$ on $G/K$, introduced in \cite{Zh}. Accordingly to \cite{Zh}, in the Siegel domain realization $T_{\Omega}$ of $G/K$ the function $e_{\lambda,\nu}$ is given by 
$$
\tilde{e}_{\lambda,\nu}(w)=\Delta_{\mu_{\lambda}+\rho}(\omega(w))\Delta(\omega(w))^{-\frac{\nu}{2}}, \quad w\in T_{\Omega}.
$$
In above $\Delta$ denotes the Koecher norm function on $\mathfrak{p}_{+}$, and $\omega(w)= \frac{w+\bar{w}}{2}$, see \cite{F} for more details.\\
Let $\gamma$ be the Cayley transform from $\mathcal{D}$ onto the Siegel realization $T_{\Omega}$ of $G/K$. Then we have 
$$
\Psi_{\lambda,\nu}(z)= \Delta(e-z)^{-\nu}\Delta_{\mu_{\lambda}+\rho}(\omega(\gamma(z)))\Delta(\omega(\gamma(z)))^{-\frac{\nu}{2}},
$$
and since 
$$
h(z,z)=\Delta(e-z)\Delta(\omega(\gamma(z))\overline{\Delta(e-z)},
$$
we get 
$$
\Psi_{\lambda,\nu}(z)=[\frac{h(z,z)}{\mid h(z,e)\mid }]^{-\frac{\nu}{2}}h(z,e)^{-\nu}.
$$
Observing that $P_{\lambda,\nu}(z,v)=\Psi_{\lambda,\nu}(k^{-1}g)$ and that $h(z,w)$ is $K$ bi-invariant we obtain the desired result.\\
ii) The identity (5.1) is easily derived from the following identity on the Jordan polynomial $h(z,w)$
$$
h(g.z,g.w)=J_{g}(z)^{\frac{1}{p}}h(z,w)\overline{J_{g}(u)^{\frac{1}{p}}},
$$ 
and the proof of Proposition 5.1 is finished.\\
It follows from above that the Poisson transform can be now rewritten explicitly as 
\begin{equation}\label{Poisson}
P_{\lambda,\nu}f(z)=\int_{S}\left(\frac{h(z,z)}{\mid h(z,u)\mid^{2}}\right)^\frac{\lambda+\eta-\nu}{2}h(z,u)^{-\nu}f(u)du,
\end{equation}

\begin{remark}
This result has been proved recently by Koranyi \cite{K2}  by using a different method.\\ Formula \ref{Poisson} agrees with Theorem 4.2 of \cite{K2}, with $l=\nu, q=\eta$ and $s=\eta-\lambda$.
\end{remark} 
Now we are ready to prove the main result of this section.

\begin{proposition}
 Let $F=P_{\lambda,\nu}f$ with  $f\in B(G/P_{\Xi},L_{\lambda,\xi})$. Then
\begin{eqnarray}
 \mathcal{H} F=\frac{(\lambda^{2}-(\eta-\nu)^{2})}{4p}F.(-iZ_{0}).
\end{eqnarray}
\end{proposition}

Proof. 
In view of Lemma 5.1 and the invariance property (5.1) of the Poisson kernel,  it suffices to show that  $\mathcal{H}_{\nu}P_{\lambda,\nu}(z,v)_{\mid z=0}=\frac{(\lambda^{2}-(\eta-\nu)^{2})}{4p}.(-iZ_{0})$.\\
To do so let $E_{j}$ be an orthonormal basis of $\mathfrak{p}_{+}$ and $E^{\ast}_{j}$ be a dual basis of $\mathfrak{p}_{-}$ with respect to $B$ (for example $\{\tilde{E}_{\alpha}\}_{\alpha\in \Phi^{+}}$ and  $\{\tilde{E}_{-\alpha}\}_{\alpha\in \Phi^{+}}$ are such basis). Let $z_{1},...,z_{n}$ be coordinates for $\mathfrak{p}_{+}$ with respect to $E_{j}$.\\

Regarding $K$-invariant functions on $G$ as functions on $\mathcal{D}$  and vice versa, we have
$$
E_{i}E^{\ast}_{j}F(e)=\frac{\partial^{2}}{\partial z_{i}\partial \bar{z_{j}}}F(0).
$$

We know that   
$$
h(z,v)^{-\nu}=1+\nu <z,v> + \mbox{higher order homogeneous terms},
$$
where $<z,v>=-\frac{1}{2p} B(z,\tau v)$.
A simple computation gives
$$
\frac{\partial^{2}}{ \partial z_{i}\partial\bar{z_{j}}}P_{\lambda,\nu}(z,v)_{\mid z=0}=\frac{(\lambda+\eta-\nu)(\lambda+\eta+\nu)}{4}v_{j}\bar{v_{i}}-\frac{\lambda+\eta-\nu}{2}\delta_{ij}.
$$
Therefore
$$
\mathcal{H}_{\nu} P_{\lambda,\nu}(z,v)_{\mid z=0}=\sum_{i,j}[E_{j},E^{\ast}_{i}](\frac{(\lambda+\eta-\nu)(\lambda+\eta+\nu)}{4}v_{j}\bar{v_{i}}-\frac{\lambda+\eta-\nu}{2}\delta_{ij})
$$
$$
=\frac{(\lambda+\eta-\nu)(\lambda+\eta+\nu)}{4}[v,\bar{v}]-\frac{\lambda+\eta-\nu}{2}\sum_{\alpha\in\Phi^{+}}\tilde{H_{\alpha}}.
$$

Notice that $[v,\bar{v}]=[Ad(k)\tilde{E_{0}},\bar{Ad(k)\tilde{{E}_{0}}}]=\sum^{r}_{j=1}\tilde{{H}_{\gamma_{j}}}$.\\
Since 
$$
\frac{1}{\eta}\sum_{\alpha\in\Phi^{+}}\tilde{H_{\alpha}}=\sum^{r}_{j=1}\tilde{{H}_{\gamma_{j}}},
$$
we get 
$$
\mathcal{H}_{\nu} P_{\lambda,\nu}(z,v)_{\mid z=0}=\frac{(\lambda+\eta-\nu)(\lambda-\eta+\nu)}{4}\tilde{H_{0}},
$$
as 
$$
\sum^{r}_{j=1}\tilde{{H}_{\gamma_{j}}}=\frac{-i}{p}Z_{0},
$$ 
the result follows.

\section{The Hua eigensections}
In this section we shall consider the  subsystem  
$$
\mathcal{H}_{\mathfrak{q}} F=\frac{(\lambda^{2}-(\eta-\nu)^{2})}{4p}F.(-iZ_{0}),
$$
and prove the main result of this section.

\begin{theorem}
Let $F\in B(G/K,\tau)$ such that 
$$
\mathcal{H}_{\mathfrak{q}} F=\frac{(\lambda^{2}-(\eta-\nu)^{2})}{4p}F.(-iZ_{0}).
$$
Then $F\in \mathcal{A}(G/K,\mathcal{M}_{\mu_{\lambda},\nu})$.
\end{theorem}
Most of the Proof of Theorem 6.1 consists in proving the following
\begin{theorem}
Let $F$ be a $\tau_{-\nu}$-spherical function on $G$ satisfying 
\begin{eqnarray}
\mathcal{H}_{\mathfrak{h}}F=\frac{\lambda^{2}-(\eta-\nu)^{2}}{4p}F.(-i)Z_{0}.
\end{eqnarray}
Then up to a constant multiple we have 
$$
F(a_{t})=\prod^{r}_{j=1}(1-\tanh^{2}t_{j})^{\frac{\lambda+\eta}{2}}
\quad_{2}F_{1}^{(m)}(\frac{\lambda+\eta-\nu}{2},\frac{\lambda+\eta+\nu}{2},\eta;\tanh^{2} t_{1},...,\tanh^{2} t_{r}),
$$
where $a_{t}=\exp(\sum^{r}_{j=1}t_{j}X_{\gamma_{j}})$.\\
In particular the elementary spherical function $\Phi_{\mu_{\lambda},\nu}$ is given by 
\begin{equation}
\Phi_{\mu_{\lambda},\nu}(a_{t})=\prod^{r}_{j=1}(1-\tanh^{2}t_{j})^{\frac{\lambda+\eta}{2}}
\quad_{2}F_{1}^{(m)}(\frac{\lambda+\eta-\nu}{2},\frac{\lambda+\eta+\nu}{2},\eta;\tanh^{2} t_{1},...,\tanh^{2} t_{r}).
\end{equation}
\end{theorem}
In the above $_{2}F_{1}^{(m)}(\alpha,\beta,\gamma;x_{1},....,x_{r})$ is the generalized Gauss hypergeometric function, see \cite{Y} for more details.\\

\begin{corollary}
let $\lambda\in{\C}$ and let $\nu\in{\Z}$. Then we have
\begin{equation}
\int_{S}[\frac{h(z,z)}{\mid h(z,u)\mid^2}]^{\frac{\lambda+\eta-\nu}{2}}h(z,u)^{-\nu}du=h(z,z)^{\frac{\lambda+\eta-\nu}{2}}
\quad_{2}F_{1}^{(m)}(\frac{\lambda+\eta-\nu}{2},\frac{\lambda+\eta+\nu}{2};\eta,\tanh^{2}t_{1},...,\tanh^{2}t_{r}),
\end{equation}
$z=Ad(k)\sum\limits_{j=1}^{r}\tanh t_{j}E_{\gamma_{j}}$.
\end{corollary}

\textbf{Proof.}
We first note that
\begin{equation}
\int_{S}[\frac{h(z,z)}{\mid h(z,u)\mid^2}]^{\frac{\lambda+\eta-\nu}{2}}h(z,u)^{-\nu}du=\tau_{\nu}(U(g:0))\Phi_{\mu_{\lambda},\nu}(g),\quad z=g.0
\end{equation}
We have
$$
\tau_{\nu}(U(g:0))\Phi_{\mu_{\lambda},\nu}(g)=\tau_{\nu}(U(a_{t}:0))\Phi_{\mu_{\lambda},\nu}(a_{t}),
$$
if $g=ha_{t}k$ with respect to the Cartan decomposition $G=KCl(A^{+})K$.\\

Now recall that the $P^{+}K_{c}P^{-}$ decomposition of $a_{t}=\exp(\sum^{r}_{j=1}t_{j}X_{\gamma_{j}})$ is
$$
a_{t}=\exp(\sum^{r}_{j=1}\tanh t_{j}E_{\gamma_{j}})\exp(\sum^{r}_{j=1}(-\log\cosh t_{j})H_{\gamma_{j}})\exp(\sum^{r}_{j=1}(\tanh t_{j}E_{-\gamma_{j}}).
$$
Hence 
$$
\tau_{\nu}(U(a_{t}:0))=(\prod\limits_{j=1}^{r}\cosh t_{j})^{\nu}.
$$
We therefore obtain
$$
\int_{S}[\frac{h(z,z)}{\mid h(z,u)\mid^2}]^{\frac{\lambda+\eta-\nu}{2}}h(z,u)^{-\nu}du=(\prod\limits_{j=1}^{r}\cosh t_{j})^{\nu}\Phi_{\mu_{\lambda},\nu}(a_{t}),
$$
and from (6.2) we obtain the identity (6.3).\\

\textbf{Proof of Theorem 6.1.} Let $F\in B(G/K,\tau_{\nu})$ such that 
\begin{eqnarray}
\mathcal{H}_{\mathfrak{q}} F=\frac{\lambda^{2}-(\eta-\nu)^{2}}{4p}F.(-iZ_{0}).
\end{eqnarray}

Fix $g\in G$ and put $F_{g}(x)=\int_{K}F(gkx)\tau_{\nu}(k)dk$. Then  $F_{g}$ satisfies  
$$
F_{g}(k_{1}xk_{2})=\tau_{\nu}^{-1}(k_{1})F(_{g}(x)\tau_{\nu}^{-1}(k_{2}),
$$

and the system (6.1). Therefore  $F_{g}=F_{g}(e)\Phi_{\mu_{\lambda},\nu}$, by Theorem 6.2.\\ 
That is  
$$
\int_{K}F(gkx)\tau_{\nu}(k)dk=\Phi_{\mu_{\lambda},\nu}(x)F(g).
$$ 
As in the case $\nu=0$ we can  prove that the above functional equation characterizes the joint eigensections of $D_{\nu}(G/K)$, so  $F\in \mathcal{A}(G/K,\mathcal{M}_{\mu_{\lambda},\nu})$ and the proof of Theorem 6.1 is finished.\\

To prove Theorem 6.2 we shall need an explicit form of the radial components of the operator $\mathcal{H}_{\mathfrak{h}}$.

\subsection{Radial components.} 
Recall that  function $F$ on $G$ is $\tau_{-\nu}$-spherical if it satisfies
$$
F(k_{1}gk_{2})=\tau_{\nu}(k_{2})^{-1}F(g)\tau_{\nu}(k_{1})^{-1}.
$$
Let $\overline{F}$ denote the restriction to $A^{+}$ of a $\tau$-spherical function $F$. By the Cartan decomposition $G=KCl(A^{+})K$ a $\tau_{-\nu}$-spherical function $F$ is essentially determined by $\overline{F}$.\\
For $U\in U(\mathfrak{g}_{\C})$, we denote by $\Delta_{\tau}(U)$ its $\tau$-radial component. That is $\Delta_{\tau}(U)\in U(\mathfrak{a}_{c})$  and satisfies
$$
\overline{(U F)}=\Delta_{\tau}(U)\overline{F},
$$
for every $\tau_{-\nu}$-spherical function $F$ on $G$.\\
Let $H^{\ast}_{\gamma_{k}}$ $(k=1,...,r)$  be the basis of $\mathfrak{h}$ which is dual to $H_{\gamma_{k}}$. Then we have
$$
\mathcal{H}_{\mathfrak{h}}=\sum^{r}_{k=1}U_{k}\otimes H^{\ast}_{\gamma_{k}},
$$
where the components $U_{k}$ are given by
$$
U_{k}=\sum_{\alpha\in \Phi^{+}}\alpha(H_{\gamma_{k}})E_{\alpha}E^{\ast}_{\alpha}.
$$
We have
 $$
 \Delta_{\tau}(U_{k})= \sum_{\alpha\in \Phi^{+}}\alpha(H_{\gamma_{k}})\Delta_{\tau}(E_{\alpha}E^{\ast}_{\alpha}).
$$
Let $F$ be a $\tau_{-\nu}$-spherical function on $G$. Then the system  $\mathcal{H}_{\mathfrak{h}}F=\frac{(\lambda^{2}-(\eta-\nu)^{2})}{4p}F.(-i)Z_{0}$ reads as
\begin{eqnarray}
(\Delta_{\tau}(U_{k}))\overline{F}(a)=\frac{(\lambda^{2}-(\eta-\nu)^{2})}{4p}F,
\end{eqnarray}
for $k=1,2,...,r$.\\
Recall that the dual basis $E^{\ast}_{\alpha}$  of $\mathfrak{p}_{-}$ is given by 
$$
E^{\ast}_{\alpha}=\frac{<\alpha,\alpha>}{2}E_{-\alpha}.
$$
In view of Proposition 2.1, $U_{k}$ can be rewritten as
$$
U_{k}=\mid\gamma_{k}\mid^{2} E_{\gamma_{k}}E_{-\gamma_{k}}+\frac{1}{2}\sum_{j}\sum_{\alpha\in \Phi^{+}_{jk}}<\alpha,\alpha>E_{\alpha}E_{-\alpha}.
$$
To compute  the $\tau$-radial parts $\Delta_{\tau}(U_{k})$ explicitly, we should determine $\Delta_{\tau}(E_{\alpha}E_{-\alpha})$ for all positive noncompact roots $\alpha $.\\

The representation $\tau_{\nu}$ of $K$ induces differentiated representation of the Lie algebra $\mathfrak{k}$. We shall denote this representation by the same latter $\tau_{\nu}$.

\begin{lemma}
Let  $F$ be a $\tau_{\nu}$-spherical function on $G$ and let $\alpha\in \Phi_{+}$. 

i) If $\alpha\in \Gamma $ or $\alpha\in \Phi^{+}_{ij}$ with $\alpha\neq \tilde{\alpha}$, then 
$$ 
\tau_{\nu}(-iH_{\alpha})F(a)=-i\nu
$$
   
ii) If $\alpha\in \Phi^{+}_{ij}$ with $\alpha=\tilde{\alpha}$, then 
$$
 \tau_{\nu}(-iH_{\alpha})F(a)=-2i\nu 
$$
 \end{lemma}

Proof. 
We have  
\begin{eqnarray}
B(-iH_{\alpha},Z_{0})=\frac{2}{\mid\alpha\mid^{2}},
\end{eqnarray}
i) If $\alpha\in \gamma_{k}$ for some $k=1,...,r$ or $\alpha\in \Phi^{+}_{ij}$ with $\alpha\neq \tilde{\alpha}$ then  
$$
 <Z-iH_{\alpha},Z_{0}>=0,
$$
by iv) of Proposition 2.1. Thus $ Z-iH_{\alpha}\in\mathfrak{t}_{s}$.
Therefore 
$$
\tau^{-1}_{\nu}(\exp t(-iH_{\alpha})=e^{-i\nu t},
$$
and i) follows.\\
ii) In the case $\alpha\in\Phi^{+}_{ij}$ with $\alpha=\tilde{\alpha}$ then 
$$
<Z-\frac{i}{2}H_{\alpha},Z_{0}>=0,
$$ 
from which we deduce that 
$
\tau^{-1}_{\nu}(\exp t(\frac{i}{2}H_{\alpha})=e^{-i\nu t}.
$
Hence ii).

\begin{proposition}
$$
4\Delta_{\tau}(E_{\gamma_{k}}E_{-\gamma_{k}})=X^{2}_{\gamma_{k}}+2\coth 2t_{k}X_{\gamma_{k}}-\nu^{2}\tanh^{2}t_{k}+2\nu.
$$
\end{proposition}

Proof.Write 
$$
E_{\gamma_{k}}E_{-\gamma_{k}}=\frac{1}{4}(X_{\gamma_{k}}^{2}+Y^{2}_{\gamma_{k}}+i[X_{\gamma_{k}},Y_{\gamma_{k}}]), 
$$
and since 
$$
[X_{\gamma_{k}},Y_{\gamma_{k}}]=-2iH_{\gamma_{k}},
$$
we get
$$
E_{\gamma_{k}}E_{-\gamma_{k}}=\frac{1}{4}(X_{\gamma_{k}}^{2}+Y^{2}_{\gamma_{k}}+i(-2iH_{\gamma_{k}})).
$$
Hence we need to compute only   $\Delta_{\tau}(Y^{2}_{\gamma_{k}})$.\\
Let $a=\exp \sum_{j=1}^{r}t_{j}X_{\gamma_{j}}$. Then we have 
$$
Ad(a^{-1})iH_{\gamma_{k}}=(\cosh 2t_{k})iH_{\gamma_{k}}+(\sinh 2t_{k})Y_{\gamma_{k}}.
$$
Thus
$$
Y^{2}_{\gamma_{k}}=(\coth 2t_{\gamma_{k}})(iH_{\gamma_{k}})^{2}+\sinh ^{-2}2t_{k}(Ad(a^{-1})iH_{\gamma_{k}})^{2}
$$
$$
-\coth 2t_{\gamma_{k}}\sinh ^{-1}2t_{k}(iH_{\gamma_{k}}Ad(a^{-1})iH_{\gamma_{k}})+(Ad(a^{-1})iH_{\gamma_{k}})iH_{\gamma_{k}}.
$$
Observe that 
$$
[Ad(a^{-1})iH_{\gamma_{k}},iH_{\gamma_{k}}]=(2\sinh 2t_{k})X_{\gamma_{k}},
$$
since $[Y_{\gamma_{k}},iH_{\gamma_{k}}]=2X_{\gamma_{k}}$.\\
Next, since $F$ is a $\tau_{\nu}$-spherical function we have 
$$
(Ad(a^{-1})iH_{\gamma_{k}}) iH_{\gamma_{k}}F(a)=-\nu^{2}F(a),
$$
hence
$$
(iH_{\gamma_{k}}Ad(a^{-1})iH_{\gamma_{k}})F(a)=-2\sinh 2t_{k}X_{\gamma_{k}}-\nu^{2},
$$
from which we deduce that 
$$
Y^{2}_{\gamma_{k}}=2\coth 2t_{k}X_{\gamma_{k}}-\nu^{2}\tanh^{2}t_{k},
$$
and the proof of Proposition 6.1 is finished.\\

\begin{proposition}
let $\alpha\in\Phi^{+}$ such that $\alpha\sim\frac{\gamma_{i}+\gamma_{j}}{2}$ ( $i\neq j$), with $\alpha\neq\widetilde{\alpha}$. Then we have
$$
 \Delta_{\tau}(E_{\alpha}E_{-\alpha}+E_{\widetilde{\alpha}}E_{\widetilde{-\alpha}})=\frac{1}{2}[\coth(t_{i}+t_{j})(X_{\gamma_{i}}+X_{\gamma_{j}})+
 \coth(t_{i}-t_{j})(X_{\gamma_{i}}-X_{\gamma_{j}})+2\nu].
 $$
\end{proposition}
We first prepare the following Lemma 
\begin{lemma}
$$
i) Ad(a^{-1})U_{1}=\cosh(t_{i}+t_{j})U_{1}-\sinh(t_{i}+t_{j})(Y_{\alpha}-\epsilon_{\alpha}Y_{\tilde{\alpha}}),
$$
$$
ii) Ad(a^{-1})U_{2}=\cosh(t_{i}+t_{j})U_{2}-\sinh(t_{i}+t_{j})(X_{\alpha}+\epsilon_{\alpha}X_{\tilde{\alpha}}),
$$
$$
iii) Ad(a^{-1})V_{1}=\cosh(t_{i}-t_{j})V_{1}-\sinh(t_{i}-t_{j})(X_{\alpha}-\epsilon_{\alpha}X_{\tilde{\alpha}}),
$$
$$
iv) Ad(a^{-1})V_{2}=\cosh(t_{i}-t_{j})V_{2}-\sinh(t_{i}-t_{j})(Y_{\alpha}+\epsilon_{\alpha}Y_{\tilde{\alpha}}).
$$
\end{lemma}
To prove the above lemma we need the following result from \cite{L} 

\begin{lemma}\cite{L}.
Let $\alpha\in \Phi_{ij}^{+}$. Then there exists $U_{1},U_{2}$ in $\mathfrak{q}$ and $V_{1},V_{2}$ in $\mathfrak{l}$ such that 
$$
i) ad(X_{\gamma_{k}})U_{1}=(\delta_{ik}+\delta_{jk})(Y_{\alpha}-\epsilon_{\alpha}Y_{\tilde{\alpha}}),
$$
$$
ii)ad(X_{\gamma_{k}})(Y_{\alpha}-\epsilon_{\alpha}Y_{\tilde{\alpha}})=(\delta_{ik}+\delta_{jk})U_{1},
$$
$$
iii) ad(X_{\gamma_{k}})U_{2}=(\delta_{ik}+\delta_{jk})(X_{\alpha}+\epsilon_{\alpha}X_{\tilde{\alpha}}),
$$
$$
iv) ad(X_{\gamma_{k}})(X_{\alpha}+\epsilon_{\alpha}X_{\tilde{\alpha}})=(\delta_{ik}+\delta_{jk})U_{2},
$$
$$
v) ad(X_{\gamma_{k}})V_{1}=(\delta_{ik}-\delta_{jk})(X_{\alpha}-\epsilon_{\alpha}X_{\tilde{\alpha}}),
$$
$$
vi) ad(X_{\gamma_{k}})(X_{\alpha}-\epsilon_{\alpha}X_{\tilde{\alpha}})=(\delta_{ik}-\delta_{jk})V_{1},
$$
$$
vii)ad(X_{\gamma_{k}})V_{2}=(\delta_{ik}-\delta_{jk})(Y_{\alpha}+\epsilon_{\alpha}Y_{\tilde{\alpha}}),
$$ 
$$
viii) ad(X_{\gamma_{k}})(Y_{\alpha}+\epsilon_{\alpha}Y_{\tilde{\alpha}})= (\delta_{ik}-\delta_{jk})V_{2}.
$$
\end{lemma}

Proof of Lemma 6.2. We shall prove only i), the others assertions can be proved in a similar way.\\
we have 
$$
Ad(\exp (-\sum^{r}_{k=1}t_{k}X_{\gamma_{k}})U_{1}=\exp(ad(-\sum^{r}_{k=1} X_{\gamma_{k}}))U_{1}.
$$
By i) of the above lemma, we have
$$
ad(-\sum^{r}_{k=1} X_{\gamma_{k}})^{2n}U_{1}=(t_{i}+t_{j})^{2n}U_{1},
$$
and 
$$
ad(-\sum^{r}_{k=1} X_{\gamma_{k}})^{2n+1}U_{1}=(t_{i}+t_{j})^{2n+1}(Y_{\alpha}-\epsilon_{\alpha}Y_{\tilde{\alpha}}),
$$
from which we get
$$
Ad(a^{-1})U_{1}=\cosh(t_{i}+t_{j})U_{1}-\sinh(t_{i}+t_{j})(Y_{\alpha}-\epsilon_{\alpha}Y_{\tilde{\alpha}}).
$$
This finishes the proof.\\
Now we  prove Proposition 6.2 giving the $\tau$-radial part of $E_{\alpha}E_{-\alpha}+E_{\widetilde{\alpha}}E_{\widetilde{-\alpha}}$.\\

Proof of Proposition 6.2. First observe that 
$$
4(E_{\alpha}E_{-\alpha}+E_{\tilde{\alpha}}E_{-\tilde{\alpha}})=X^{2}_{\alpha}+X^{2}_{\tilde{\alpha}}+
Y^{2}_{\alpha}+Y^{2}_{\tilde{\alpha}}+i(2iH_{\alpha}+2iH_{\tilde{\alpha}})).
$$
Next, since
$$X^{2}_{\alpha}+X^{2}_{\tilde{\alpha}}=\frac{1}{2}[(X_{\alpha}+\epsilon_{\alpha} X_{\tilde{\alpha}})^{2}+(X_{\alpha}-\epsilon_{\alpha} X_{\tilde{\alpha}})^{2}],
$$
we will compute $\Delta_{\tau}((X_{\alpha}\pm\epsilon_{\alpha} X_{\tilde{\alpha}})^{2}$.\\
To this end, consider the element $U_{2}$ given by Lemma 6.3 and let $a=\exp(\sum^{r}_{k=1}t_{j}X_{\gamma_{j}})$. Then we have
$$
Ad(a^{-1})U_{2}=\cosh(t_{i}+t_{j})U_{2}-\sinh(t_{i}+t_{j})(X_{\alpha}+\epsilon X_{\tilde{\alpha}}),
$$
from which we get 
$$
\sinh^{2}(t_{i}+t_{j})(X_{\alpha}+\epsilon X_{\tilde{\alpha}})^{2}=\cosh^{2}(t_{i}+t_{j})U^{2}_{2}-\cosh(t_{i}+t_{j})U_{2}Ad(a^{-1})U_{2}-
$$
$$
\cosh(t_{i}+t_{j})(Ad(a^{-1})U_{2})U_{2}+(Ad(a^{-1})U_{2})^{2}.
$$
Recall from \cite{L}, that $U_{2}=i(Q_{\alpha}-\bar{Q_{\alpha}})$ where
$$
Q_{\alpha}=N_{\alpha,-\gamma_{i}}E_{\alpha-\gamma_{i}}+N_{\alpha,-\gamma_{j}}E_{\alpha-\gamma_{j}},
$$
from which we obtain  $<U_{2},Z_{0}>=0$, therefore $U_{2}\in \mathfrak{q}\cap \mathfrak{t}_{s}$ and $\tau_{\nu}(U_{2})=0$ .\\
By ii) in Lemma 6.2, we have
$$
[Ad(a^{-1})U_{2}),U_{2}]=-2 \sinh(t_{i}+t_{j})(X_{\gamma_{i}}+X_{\gamma_{j}}),
$$
which gives us 
$$
\Delta_{\tau}((X_{\alpha}+\epsilon_{\alpha} X_{\tilde{\alpha}})^{2}=2\coth(t_{i}+t_{j})(X_{\gamma_{i}}+X_{\gamma_{j}}).
$$
Similarly by considering $V_{1}$ and noticing that $V_{1}\in\mathfrak{l}\cap \mathfrak{t}_{s}$ we get  
$$
\Delta_{\tau}((X_{\alpha}-\epsilon_{\alpha} X_{\tilde{\alpha}})^{2}=2\coth(t_{i}-t_{j})(X_{\gamma_{i}}-X_{\gamma_{j}}).
$$
Finally since $-2iH_{\tilde{\alpha}}F(a)=-2i\nu F(a)$,by Lemma 6.1, the result follows. 

\begin{proposition}
Let $\alpha\in\Phi^{+}$ such that $\alpha=\frac{\gamma_{i}+\gamma_{j}}{2}$ $(i\neq j)$, that is $\alpha=\widetilde{\alpha}$. Then we have
$$
\Delta_{\tau}(E_{\alpha}E_{-\alpha})=\frac{1}{2}[\coth(t_{i}+t_{j})(X_{\gamma_{i}}+X_{\gamma_{j}})+
 \coth(t_{i}-t_{j})(X_{\gamma_{i}}-X_{\gamma_{j}})+2\nu].
$$
\end{proposition}

Proof. We first suppose that $\epsilon_{\alpha}=1$. Accordingly to Proposition 13 and Proposition 14 in \cite{L}, $U_{1}=V_{1}=0$.\\
We compute first $\Delta_{\tau}(X_{\alpha}^{2})$ and $\Delta_{\tau}(Y_{\alpha}^{2})$ in the case $\epsilon_{\alpha}=1$.\\
We have
$$
Ad(a^{-1})U_{2}=\cosh(t_{i}+t_{j})U_{2}-2\sinh(t_{i}+t_{j})X_{\alpha},
$$
hence
$$
4\sinh^{2}(t_{i}+t_{j})X^{2}_{\alpha}=\cosh^{2}(t_{i}+t_{j})U^{2}_{2}+(Ad(a^{-1})U_{2})^{2}-\cosh(t_{i}+t_{j})(U_{2}Ad(a^{-1})U_{2}+(Ad(a^{-1})U_{2})U_{2}).
$$
Noticing that 
$$
U_{2}Ad(a^{-1})U_{2}=-[Ad(a^{-1})U_{2},U_{2}]
$$
$$
=-8\sinh(t_{i}+t_{j})(X_{\gamma_{i}}+X_{\gamma_{j}}),
$$
we get
$$
\Delta_{\tau}(X_{\alpha}^{2})=2\coth(t_{i}+t_{j})(X_{\gamma_{i}}+X_{\gamma_{j}}),
$$
and similarly 
$$
\Delta_{\tau}(Y_{\alpha}^{2})=2\coth(t_{i}-t_{j})(X_{\gamma_{i}}-X_{\gamma_{j}}).
$$
in the case $\epsilon_{\alpha}=-1$, analogous computations give
$$
\Delta_{\tau}(Y_{\alpha}^{2})=2\coth(t_{i}+t_{j})(X_{\gamma_{i}}+X_{\gamma_{j}}),
$$
and 
$$  
\Delta_{\tau}(X_{\alpha}^{2})=2\coth(t_{i}-t_{j})(X_{\gamma_{i}}-X_{\gamma_{j}}).
$$
To finish the proof of Proposition 6.3, notice that in the case $\alpha= \tilde{\alpha}$
$$
-i2H_{\alpha}F(a)=-4i\nu.
$$

\begin{proposition}
The $\tau$-radial part of the operator $U_{k}$ is given by
$$
\frac{4}{\mid\gamma_{k}\mid^{2}}\Delta_{\tau}(U_{k})=\frac{\partial^{2}}{\partial t^{2}_{k}}+2\coth2t_{k}\frac{\partial}{\partial t_{k}}-\nu^{2}\tanh^{2}t_{k}+2\nu
$$
$$
+\frac{m}{2}\sum^{r}_{j\neq k}[\coth(t_{j}+t_{k})(\frac{\partial}{\partial t_{j}}+\frac{\partial}{\partial t_{k}})+
\coth(t_{j}-t_{k})(\frac{\partial}{\partial t_{j}}-\frac{\partial}{\partial t_{k}})+2\nu].
$$
\end{proposition}

Proof.
Recall that
$$
\Delta_{\tau}(U_{k})= \sum_{\alpha\in \Phi^{+}}\alpha(H_{\gamma_{k}})\Delta_{\tau}(E_{\alpha}E^{\ast}_{\alpha}).
$$
By using (i), (ii) and (iii) of proposition 2.1, we obtain
$$
\Delta_{\tau}(U_{k})=<\gamma_{k},\gamma_{k}>\Delta_{\tau}(E_{\gamma_{k}}E_{-\gamma_{k}})+\frac{1}{2}\sum_{j\neq k}\sum_{\alpha\in\Phi_{jk}^{+}}<\alpha,\alpha>
\Delta_{\tau}(E_{\alpha}E_{-\alpha}),
$$
from which we deduce  
$$
\frac{4}{\mid\gamma_{k}\mid^{2}}\Delta_{\tau}(U_{k})=\Delta_{\tau}(E_{\gamma_{k}}E_{-\gamma_{k}})+\frac{m}{2}\sum_{j\neq k}
[\coth(t_{j}+t_{k})(X_{\gamma_{j}}+X_{\gamma_{k}})+
 \coth(t_{i}-t_{j})(X_{\gamma_{j}}-X_{\gamma_{k}})+2\nu],
$$
by the results of Proposition 6.1, Proposition 6.2 and Proposition 6.3. 

Next consider  a coordinate system $t=(t_{1},...,t_{r})\in{\R}^{r} \longmapsto \exp(\sum^{r}_{j=1}t_{j}X_{\gamma_{j}})\in A $, such that the Weyl group $W$ acts as the group of all permutations and sign changes of the coordinates $(t_{1},...,t_{r})$ to get the result. This finishes the proof of Proposition 6.4. \\
\subsection{Proof of Theorem 6.2}
In this subsection we give the proof Theorem 6.2 which is the main step in the proof of our main result. 

Proof of Theorem 6.2 . Let $F$ be a $\tau_{-\nu}$-spherical function satisfying the system (6.1). \\
Recall from subsection 6.1 (equation (6.6) ) that $F$ satisfies
$$
(\Delta_{\tau}(U_{k}))\overline{F}(a)=\frac{(\lambda^{2}-(\eta-\nu)^{2})}{4p}F, \quad k=1,...,r
$$
Let 
$$
\phi(t_{1},...,t_{r})=(\prod^{r}_{j=1}\cosh t_{j})^{\nu}F(t_{1},...,t_{r}).
$$
It follows from Proposition 6.4 that  the function $\phi(t_{1},...,t_{r})$ satisfies the system of differential equations 
$$
\frac{\partial^{2}\phi}{\partial t^{2}_{k}}+2\coth2t_{k}\frac{\partial\phi}{\partial t_{k}}-2\nu\tanh t_{k}\frac{\partial\phi}{\partial t_{k}}
+\frac{m}{2}\sum^{r}_{j=1}\frac{1}{(\sinh^{2}t_{j}-\sinh^{2}t_{k})}(\sinh 2t_{j}\frac{\partial\phi}{\partial t_{j}}-\sinh 2t_{k}\frac{\partial\phi}{\partial t_{k}})
$$
$$
=\frac{(\lambda^{2}-(\eta-\nu)^{2})}{4}\phi,
$$
for all $k=1,...,r$.\\

Put $x_{i}=-\sinh^{2}t_{i}$ and  $\psi(x_{1},....,x_{r})=\phi((t_{1},...,t_{r})$. Then the function $\psi$ satisfies the system
 $$
 x_{k}(1-x_{k})\frac{\partial^{2}\psi}{\partial x^{2}_{k}}+(1-(2-\nu)x_{k})\frac{\partial \psi}{\partial x_{k}}
 -\frac{m}{2}\sum_{j\neq k}\frac{x_{j}(1-x_{j})}{x_{k}-x_{j}}\frac{\partial \psi}{\partial x_{j}}-\frac{x_{k}(1-x_{k})}{x_{k}-x_{j}}
 \frac{\partial \psi}{\partial x_{k}}
 $$
 $$
 =\frac{(\eta-\nu)^{2}-\lambda^{2}}{4}\psi,
$$
for all $k=1,...,r$.\\
Since $\mathfrak{z}\subset \mathfrak{h}$, $F$ is an eigenfunction of the Laplace-Beltrami operator (see the proof of Theorem 6.1). Hence $F$ is analytic.
Being $\tau$-spherical $F$ is $W$-invariant. Since the Weyl group $W$ acts as the group of all permutations and sign changes of the coordinates $(t_{1},...,t_{r})$, it follows that $\psi$  is a symmetric function of $x_{1},...,x_{r}$ and analytic at $x_{1}=...=x_{r}=0$.\\ From Theorem 2.1 \cite{Y} we deduce that
$$
\psi(x_{1},...,x_{r})=c \quad_{2}F_{1}^{(m)}(\frac{\lambda+\eta-\nu}{2},\frac{-\lambda+\eta-\nu}{2},\eta; x_{1},...,x_{r}),
$$
where $c$ is some numerical constant.\\
Thus 
$$
F(a_{t})= (\prod^{r}_{j=1}\cosh t_{j})^{-\nu}\quad_{2}F_{1}^{m}(\frac{\lambda+\eta-\nu}{2},\frac{-\lambda+\eta-\nu}{2},\eta;-\sinh^{2} t_{1},...,-\sinh^{2} t_{r}).
$$
 Next using the following Formula  on the generalized Gauss hypergeometric function  \cite{Y}
$$
_{2}F_{1}^{(m)}(\alpha,\beta,\gamma;y_{1},....,y_{r})=\prod^{r}_{j=1}(1-y_{j})^{-a}\quad_{2}F_{1}^{(m)}(\alpha,\gamma-\beta,\gamma;\frac{y_{1}}{y_{1}-1},....,\frac{y_{r}}{y_{r}-1}).
$$
we get
$$
F(a_{t})=\prod^{r}_{j=1}(1-\tanh^{2}t_{j})^{\frac{\lambda+\eta}{2}}
\quad_{2}F_{1}^{(m)}(\frac{\lambda+\eta-\nu}{2},\frac{\lambda+\eta+\nu}{2},\eta;\tanh^{2} t_{1},...,\tanh^{2} t_{r}),
$$
and the proof of Theorem 6.2 is finished.\\

\section{The sufficiency of the conditions}
In this section we shall complete the proof of our main result.

Let $F\in \mathcal{B}(G/K,\tau_{\nu})$ such that $F$ satisfies the Hua system
\begin{eqnarray}
\mathcal{H}F=\frac{(\lambda^{2}-(\eta-\nu)^{2})}{4p}F.(-i)Z_{0}.
\end{eqnarray}
Then $F\in \mathcal{A}(G/K,\mathcal{M}_{\mu_{\lambda},\nu})$, by Theorem 6.1.\\
By Theorem 3.1, it suffices to prove that the boundary value $\tilde{\beta}_{\mu_{\lambda},\nu}F$ is in $\mathcal{B}(G/P_{\Xi},L_{\lambda,\nu})$. To do so, we will show  that the induced equations of the subsystem $\mathcal{H}_{s}F=0$ of (7.1)for boundary values on $G/P$ characterize $\mathcal{B}(G/P_{\Xi},L_{\lambda,\nu})$, see \cite{K} for more details on induced equations that boundary values satisfy.\\
The method of the proof for $\nu=0$ in \cite{Sh} can be generalized to our situation, see also \cite{L}.\\ Below, we give an outline of the proof.  Let $\mathfrak{s}$ be the orthogonal complement of $\mathfrak{h}$ in $\mathfrak{q}$  with respect to $B$. Denote by $C^{+}$ the set of positive compact roots $\beta$ such that $\beta\sim \frac{\gamma_{j}-\gamma_{i}}{2}$ for some $j>i$.\\Let  $\{S^{\ast}_{\beta},\beta \in C^{+}\}$ be the  basis of $\mathfrak{s}_{c}$ that is dual to $\{S_{\beta},\beta \in C^{+}\}$, with respect to $B$. 
We have
$$
\mathcal{H}_{\mathfrak{s}}=\sum_{\beta\in C^{+}} U_{\beta}\otimes S^{\ast}_{\beta},
$$
where $U_{\beta}\in \mathcal{U}(\mathfrak{g}_{\C})$ is given by 
$$
U_{\beta}=\sum_{\alpha\in\Phi^{+}}[S_{\beta},E_{\alpha}]E^{\ast}_{\alpha}.
$$
The condition $\mathcal{H}_{\mathfrak{s}}F=0$ implies

$$
U_{\beta}F=0, \quad \forall \beta\in C^{+},
$$

Consider the Poincar\'e-Birkhoff-Witt Theorem  decomposition
\begin{eqnarray}
U(\mathfrak{g}_{\C})=U(\mathfrak{n}^{-}_{c}+\mathfrak{a}_{c})+\sum_{X\in \mathfrak{k}_{\C}}U(\mathfrak{g}_{\C})(X-\tau_{\nu}(X)),
\end{eqnarray}
and let $\Pi_{1}$ denote the projection of $U(\mathfrak{g}_{\C})$  to $U(\mathfrak{n}^{-}_{c}+\mathfrak{a}_{c})$ with respect to the decomposition (7.2).\\
For $U\in U(\mathfrak{g}_{\C})$, let $\tilde{U}$ be the differential operator on $N^{-}\times {\R}^{r}$ defined by $\tilde{U}=P(\Pi_{1}(U))$ where $P$ denotes the operator defined by $P(X_{\alpha})=t^{\alpha}X_{-\alpha}$ and $P(H_{j})=-t_{j}\frac{\partial}{\partial t_{j}}$, under the isomorphism from ${\R}^{r}$ onto $A$ given by
$t\rightarrow a(t)=\exp(-\sum_{t_{j}\neq 0} \log \mid t_{j}\mid H_{j})$.\\
By the Iwasawa decomposition $G=N^{-}AK$, the restriction map from $G$ to $N^{-}A$ gives an isomorphism from $\mathcal{B}(G/K,\tau_{\nu})$ to $\mathcal{B}(N^{-}A)$. \\
Now  $F$ regarded as an element of $\mathcal{B}(N^{-}A)$ satisfies the differential equation
$$
\tilde{U}_{\beta}F=0, \quad \beta\in C^{+}.
$$
Fix $\beta\in C^{+}$ such that $\beta \sim \frac{\gamma_{i}-\gamma_{i-1}}{2}$ \quad $(2\le i\le r)$.\\

Similar computations as in [\cite{Sh}, Proposition 4.4 and Proposition 4.5] show that the operator 
\begin{eqnarray}
t^{\frac{-1}{2}(\beta_{i}-\beta_{i-1})}\tilde{U}_{\beta},
\end{eqnarray}
is well defined on $N^{-}\times {\R}$ and has analytic coefficients near $t=0$,and that 
the induced equations for the system $\mathcal{H}_{\mathfrak{s}}F=0$ are 
$$
Ad c(E_{-\tilde{\beta}})\mathcal{\beta}_{\mu_{\lambda},\nu}(F)=0, \quad \forall \beta \in C^{+}; \beta \sim \frac{\gamma_{i}-\gamma_{i-1}}{2}  \quad (2\le i\le r).
$$
To conclude recall that the vectors $\{ Ad c(E_{-\tilde{\beta}}),\beta \sim \frac{\gamma_{i}-\gamma_{i-1}}{2}\}$ span the root space $g_{\frac{\beta_{i}-\beta_{i-1}}{2}}$ and that  
$\{ \frac{\beta_{i}-\beta_{i-1}}{2},2\le i\le r \}$ are the simple roots of  $ \{ \frac{\beta_{i}-\beta_{j}}{2}, 1\le j<i\le r$.  
Thus
$$
X_{\alpha} \mathcal{\beta}_{\mu_{\lambda},\nu}(F)=0 \quad \forall \alpha; \quad \alpha=\frac{\beta_{i}-\beta_{j}}{2},
$$
with $1\le j<i\le r$.\\
This shows that $\tilde{\beta}_{\mu_{\lambda},\nu}F\in \mathcal{B}(G/P_{\Xi},L_{\lambda,\nu})$ and the proof the main result of this paper is finished.

\end{document}